\documentclass[smallextended,envcountsame,envcountsect]{svjour3}
\smartqed  

\usepackage{amsfonts,dsfont,scalerel,braket,booktabs}

\usepackage[usenames, dvipsnames]{color}
\definecolor{rred}{rgb}{0.7,0,0.1}
\newcommand{\mk}{\color{black}}
\newcommand{\mkr}{\color{black}}
\definecolor{greenrb}{rgb}{0.2,0.6,0.2}
 
\usepackage{newtxtext}

\spnewtheorem*{xproof}{}{\itshape}{\rmfamily}

\renewenvironment{proof}[1][\proofname]
 {\renewcommand\xproofname{#1}\xproof}
 {\endxproof}

\usepackage{amssymb,amscd}

\usepackage{amssymb}
\usepackage{subcaption}
\usepackage{bbm}
\usepackage{empheq}
\usepackage{fancyhdr}
\usepackage{tikz}
\usetikzlibrary{arrows,positioning} 
\usepackage{dsfont}
\usepackage{needspace}
\usepackage{booktabs}
\usepackage{framed}
\usepackage{sidecap}
\usepackage{bm} 
    \usepackage{cases} 
    \usepackage{color}
\usepackage{graphicx, graphics}
\usepackage[ruled]{algorithm2e}

\usepackage{lineno}
\usepackage[breaklinks=true,colorlinks=true,linkcolor=blue,urlcolor=blue,citecolor=blue]{hyperref}

\newtheorem{thm}{Theorem}
\newtheorem{lem}{Lemma}
\newtheorem{defi}{Definition}

\newtheorem{cor}{Corollary}
\newtheorem{rem}{Remark}

\newcommand\at[2]{\left.#1\right|_{#2}}

\def\bt{\begin{thm}}
\def\et{\end{thm}}
\def\bl{\begin{lem}}
\def\el{\end{lem}}
\def\bd{\begin{defi}}
\def\ed{\end{defi}}
\def\bc{\begin{cor}}
\def\ec{\end{cor}}
\def\bp{\begin{proof}}
\def\ep{\end{proof}}
\def\br{\begin{rem}}
\def\er{\end{rem}}

\def\bpp{\begin{proof}}
\def\epp{\end{proof}}
\numberwithin{equation}{section} 
\def\be{\begin{equation}}
\def\ee{\end{equation}}
\def\bes{\begin{equation*}}
\def\ees{\end{equation*}}
\def\bea{\begin{equation} \begin{aligned}}
\def\eea{\end{aligned} \end{equation}}
\def\beas{\begin{equation*} \begin{aligned}}
\def\eeas{\end{aligned} \end{equation*}}
\def\d{\, \mathrm{d}}

\def\cX{\mathcal X}
\def\m{\mathfrak{m}}

\def\V{V}

\def\Re{{\rm Re \,}}
\def\Im{{\rm Im \,}}

\def\G{\mathcal{G}}
\def\cL{\mathcal L}

\def\bi{\begin{itemize}}
\def\ei{\end{itemize}}

\begin{document}

\title{Ruelle-Pollicott Resonances of Stochastic Systems in Reduced State Space. Part I: Theory}

\author{Micka\"el D. Chekroun \and Alexis Tantet   \and Henk A. Dijkstra \and J. David Neelin}
\institute{M.D. Chekroun \at
 Department of Earth and Planetary Sciences, Weizmann Institute, Rehovot 76100, Israel;
  Department of Atmospheric and Oceanic Sciences, University of California, Los Angeles, CA 90095-1565, USA\\
    \email{mchekroun@atmos.ucla.edu}
    \and
A. Tantet \at \'Ecole Polytechnique, Laboratoire de M\'et\'eorologie Dynamique, Palaiseau, France\\
    \and
  H.A. Dijkstra \at Institute of Marine and Atmospheric Research, Department of Physics and Astronomy,  University of Utrecht, Utrecht, The Netherlands
\and
  J.D. Neelin \at
  Department of Atmospheric and Oceanic Sciences, University of California, Los Angeles, CA 90095-1565, USA}
\date{\today}

\maketitle

\begin{abstract}
A theory of Ruelle-Pollicott (RP) resonances for stochastic differential systems is presented.  These resonances are defined as the eigenvalues of the generator (Kolmogorov operator)
of a given stochastic system.  By relying on the theory of Markov semigroups, decomposition formulas of correlation functions and power spectral densities (PSDs) in terms of RP resonances are then derived. 
These formulas describe, for a broad class of stochastic differential equations (SDEs), how the RP resonances characterize the decay of correlations as well as the signal's oscillatory components manifested by peaks in the PSD.    

It is then shown that a notion reduced RP resonances can be rigorously defined, as soon as the dynamics is partially observed within a reduced state space $V$.     
These reduced resonances are obtained from the spectral elements of reduced Markov operators acting on functions of the state space $V$, and can be estimated from series. They inform us about the spectral elements of some coarse-grained version of the SDE generator.  When the time-lag at 
which the transitions are collected from partial observations in $V$, is either sufficiently small or large, it is shown that the reduced RP resonances approximate  the (weak) RP resonances of the generator of the conditional expectation in $V$, i.e.~the optimal reduced system in $V$ obtained by averaging out the contribution of the unobserved variables. The approach is illustrated on a stochastic slow-fast system for which it is shown that the reduced RP resonances allow for a good reconstruction of the correlation functions and PSDs, even when the time-scale separation is weak. 

The companions articles, Part II \cite{PartII} and Part III \cite{PartIII}, deal with further practical aspects of the theory presented in this contribution. 
One important byproduct consists of the diagnosis usefulness of stochastic dynamics that  RP resonances provide. This is illustrated in the case of a stochastic Hopf bifurcation   
in Part II. There, it is shown that such a bifurcation has a clear manifestation in terms of a geometric organization of the RP resonances along discrete parabolas in the left half plane. 
Such geometric features formed by (reduced) RP resonances are extractable from time series and allow thus for providing an unambiguous ``signature'' of nonlinear oscillations embedded within a stochastic background.
By relying then on the theory of reduced RP resonances presented in this contribution,  Part III   addresses the question of detection and 
characterization of such oscillations in a high-dimensional stochastic system,  namely the Cane-Zebiak model of El Ni\~no-Southern Oscillation subject to noise modeling fast atmospheric fluctuations. 
\end{abstract}


\tableofcontents

\section{Introduction}

The determination of reduced equations that aims at mimicking in a reduced state space the dynamics (reduced systems for short) issued from a system of {\it stochastic differential equations (SDEs)} posed in a higher-dimensional state space has become a central question in several disciplines ranging from molecular dynamics \cite{schutte1999direct,schutte1999conformational,schutte2001transfer,zhang2017effective,duong2018quantification} to epidemic models \cite{FS09,forgoston2009accurate} and climate dynamics over the past decades; see e.g.~\cite{penland1995optimal,MajdaEtAl_CPAM01,givon2004extracting,CLW15_vol2,gottwald_crommelin_franzke_2017,berner2017stochastic} and references therein. Various approaches have been proposed, and in most of the cases, different metrics associated with moment estimates, probability density functions or time marginals are used to compare the reduced dynamics obtained from the surrogate system to that of the original system of stochastic differential equations (SDEs); see e.g.~\cite{legoll2010effective,legoll2017pathwise,zhang2017effective,duong2018quantification,CLM19_closure}. 

Here, we adopt a completely different approach which consists of using correlation functions and power spectral densities as ``metrics'' of comparison. 
In that respect, we frame hereafter a theory of Ruelle-Pollicott (RP) resonances for stochastic differential systems. Initially introduced for discrete and continuous chaotic deterministic systems \cite{ruelle1986locating,pollicott1986meromorphic}, these resonances are known to provide rigorous decomposition formulas of correlation functions and power spectra; see e.g.~\cite{baladi1989resonances,froyland1997computer,Baladi2000,melbourne2007power,giulietti2013anosov,nonnenmacher2015decay} and references therein.
However, the rigorous derivation of such formulas in the deterministic context is made difficult by the singular character of the underlying invariant measures. In presence of noise, smoothing effects are known to occur for a broad class of SDEs (see Appendix \ref{Sec_Feller}), and the derivation of such formulas is thus facilitated. Section \ref{sec:ErgodicSpectrumSDE} presents an extension of such decomposition formulas to the stochastic context.

The RP resonances are then defined naturally as the eigenvalues of the generator (Kolmogorov operator) of a given stochastic system. As shown in Sec.~\ref{sec:ErgodicSpectrumSDE}, these resonances allow for decomposition formulas of correlation functions and power spectral densities (PSDs) in terms of these resonances; see Corollary \ref{Cor_decomp} and formula \eqref{eq:spectralPower} below. The  resulting formulas describe how the RP resonances characterize the decay of correlations as well as the signal's oscillatory components manifested by peaks in the PSD.    The obtention of these formulas relies on tools and concepts from stochastic analysis, based on the theory of Markov semigroups on one hand, and the spectral theory of semigroups, on the other. These tools are briefly surveyed in Appendix \ref{Sec_Appendix}, having in mind a wider audience in the geosciences and macroscopic physics.

From a practical viewpoint, our decomposition formulas benefit furthermore from a natural dynamical interpretation as they relate to the spectral elements of the Kolmogorov operator and thus to the SDE's coefficients themselves. Such dynamical interpretations (and thus decompositions) are potentially useful for identifying physical processes responsible for power excess bumps or other broad band peaks in the PSD of noisy observations, a topic of active research in various areas of physics such as asteroseismology \cite{kallinger2010asteroseismology,gruberbauer2009detection}, supersonic flows \cite{tam1995supersonic},  or climate dynamics \cite{ghil40advanced}.

In Section \ref{Sec_Theory_reduced_spec}, we present the main contribution of this article, by inquiring whether one can extract useful resonances from 
partial observations in a reduced state space, that still relate to the correlation functions and PSDs.  
First, we show that a notion reduced RP resonances can be indeed rigorously framed, as soon as the (stochastic) dynamics is partially observed within a reduced state space $V$.   
These reduced resonances are obtained from the spectral elements of reduced Markov operators acting on functions of the state space $V$, and can be estimated from series; see Sec.~\ref{Sec_reducedRP}. 
The reduced Markov operators are ensured by Theorem \ref{Theorem_A_stoch} that extends  \cite[Theorem A]{Chek_al14_RP} within the stochastic context of this article. 
 {\it Mutatis mutandis},  Theorem \ref{Theorem_A_stoch} as \cite[Theorem A]{Chek_al14_RP}, shows that --- given a reduced state space $\V$, a mapping $h$ from the full state $\cX$ space  onto $\V$, and a Markov semigroup $P_t$ possessing an invariant measure on $\cX$ --- a family of reduced Markov operators (acting on functions of $\V$) can be naturally associated with $P_t$, $\V$ and $h$. This family characterizes a coarse-graining in the reduced state space $\V$ --- such as induced by the map $h$ --- of the actual transition probabilities associated with $P_t$ that take place within the full state space $\cX$. 
 Theorem \ref{Theorem_A_stoch} can be viewed as providing the analogue, in the SDE context dealt with in this article, of coarse-grained Markov representations used  for describing coarse-grained dynamics of Hamiltonian systems \cite{schutte1999direct} or more general Markov state models \cite{schutte2013metastability}. {\mkr We emphasize that, in the context of homogenization of SDEs \cite{pavliotis2008multiscale}, a notion of reduced RP resonances was also considered in  \cite{crommelin2011diffusion} as eigenvalues of a coarse-grained  diffusion processes obtained from multiscale data.}

 In a second step, under ergodicity (or mixing) assumptions satisfied for a broad class of SDEs, we show that the reduced RP resonances relate to the RP resonances of the generator of the conditional expectation in $V$, i.e.~the optimal reduced system in $V$ obtained by averaging out the contribution of the unobserved variables; see Theorem \ref{Thm_assymp_gene} and {\mkr Corollary \ref{Cor_approx_RP}}.  Thus, reduced RP resonances  inform us about the spectral elements of some coarse-grained version of the SDE generator. {\mkr By relying on Theorems \ref{Theorem_A_stoch} and \ref{Thm_assymp_gene}, these  resonances can be estimated from time series by estimation of transition probabilities at a time-lag $\tau$}; see Sec.~\ref{Sec_reducedRP}.  {\mkr Corollary  \ref{Cor_approx_RP} shows that an approximation relationship between the reduced RP resonances and those of the conditional expectation generator} is ensured when the time-lag $\tau$ at 
which the transitions are collected  from partial observations in $V$, is either sufficiently small or large.\footnote{In practice however it is often observed that the reduced RP resonances still provide useful information for ``intermediate" time-lags; see Part III \cite{PartIII}.}
 In the context of  homogenization of SDEs, the reduced RP resonances approximate even the RP resonances of the full SDE system; see \cite{crommelin2011diffusion} and Remark \ref{Rmk_approx_RPs}-(i).
The facts that  the reduced RP resonances can be estimated from time series and relate to the conditional expectation, provide useful insights for the investigation of whether non-Markovian effects {\it \`a la} Mori-Zwanzig should be included or not for the derivation of efficient reduced systems; see e.g.~\cite{givon2004extracting,Chorin_Hald-book,wouters2012,wouters2013,KCG_MSM15,CLM16_Lorenz9D,CLM19_closure}.      
 
In that perspective, we analyze in Sec.~\ref{Sec_Langevin} correlation functions and PSDs as reconstructed from reduced RP resonances of a stochastic slow-fast system, by using the formulas of Sec.~\ref{sec:ErgodicSpectrumSDE} in which the reduced RP resonances replace the genuine eponymous resonances. In the case of a strong time-scale separation, without any surprise these reduced RP resonances, as associated with the conditional expectation --- well-approximated by a slow manifold reduction for the example of Sec.~\ref{Sec_Langevin} ---, provide excellent reproduction of correlation functions  and PSDs obtained by standard sample estimates.  The surprise arises when the system is placed in a regime without time-scale separation. In this case, the  reduced RP resonances allow still for reconstructing to a very good accuracy level the  correlation functions  and PSDs. Thus here the conditional expectation is sufficient to close the system (at least for reconstructing correlations) and the inclusion of non-Markovian effects for model reduction can be neglected.  
Such a diagnosis regarding the conditional expectation drawn from RP resonance analysis offers thus promising perspectives in terms of reduced-order modeling, and provide useful insights to reduction approaches of SDEs exploiting conditional expectations such as in e.g.~\cite{weinan2004metastability,legoll2010effective,lu2014exact,zhang2017effective,duong2018quantification}.

The usefulness of  reduced RP resonance analysis is not limited to providing useful diagnosis for model reduction. 
The companions articles, Part II \cite{PartII} and Part III \cite{PartIII}, deal with further practical aspects of the theory presented in this article. 
One important byproduct consists of the diagnosis usefulness of stochastic dynamics that  RP resonances themselves, offer. This is illustrated in the case of a stochastic Hopf bifurcation   
in Part II. There, it is shown that such a bifurcation has a clear ``signature'' in  terms of the geometric organization of the RP resonances in the left half plane, these resonances forming in particular discrete parabolas whose characteristics describe the frequency of the underlying oscillations and whose associated modes inform about the interactions between the noise and nonlinear effects; see \cite{PartII}.

As pointed out in Sec.~\ref{Sec_reducedRP} and substantiated in Part III, 
such features are furthermore useful to identify and qualify from time series, the presence of nonlinear oscillations embedded within a stochastic background. Indeed by relying on the theory of reduced RP resonances presented in this contribution,  Part III addresses the question of detection and 
characterization of such oscillations in a high-dimensional stochastic system,  namely a system obtained from a semispectral approximation of the Cane-Zebiak model of El Ni\~no-Southern Oscillation\footnote{We refer to \cite{cao2019mathematical} for a mathematical analysis of the related Jin-Neelin model.},  with the addition of noise modeling fast atmospheric fluctuations. Finally, we mention that the characterization of nonlinear oscillations in terms of (reduced) RP resonances has already demonstrated its usefulness for data-driven modeling purposes of multiscale datasets arising in climate dynamics. Indeed by leveraging on the understanding gained here and in Part II \cite{PartII} as well as by exploiting the  framework of \cite{CK17} to decompose the multiscale variability,  successful applications to the emulation of wind-driven ocean gyres dynamics have been reported in \cite{KCB18} by means of frequency-ranked networks of stochastic Stuart-Landau oscillators. The modeling of Arctic sea ice concentration and extent has been also addressed within this framework; see \cite{KCYG17,MASIE_paper}. We envision a further program of using reduced RP resonances for the stochastic inverse modeling of climate datasets along with the analysis of the nearing to tipping points \cite{tantet2018crisis} from such datasets, in the context of climate change.

\section{Ruelle-Pollicott resonances and decay of correlations from stochastic differential equations}
\label{sec:ErgodicSpectrumSDE}

The authors in~\cite{Chek_al14_RP} have introduced a new mathematical framework to (i) understand and diagnose --- through partial observations --- the variability of chaotic flows, and (ii) to analyze parameter sensitivity that may occur in the modeling of such observations.  The framework relied on the theory of Ruelle-Pollicott (RP) resonances introduced in the mid-80's~\cite{ruelle1986locating,pollicott1986meromorphic} and that, at the time of the publication of~\cite{Chek_al14_RP}, was known only by a group of experts  working in the field of dynamical system theory \cite{Baladi2000} and the mathematical study of scattering resonances~\cite{Zworski2017}.  Initially introduced for discrete and continuous chaotic deterministic systems (see also \cite{giulietti2013anosov} for the case of Anosov flows), the RP resonances extend to stochastic differential systems. This section presents such an extension by 
relying on tools from stochastic analysis, based on the theory of Markov semigroups on one hand, and the spectral theory of semigroups, on the other; see also \cite{gaspard2002trace,dyatlov2015stochastic} for complementary approaches.

As we will see, the RP resonances characterize important features of the solution's variability (e.g.~its oscillatory behavior), such as typically reflected in power spectra or correlation functions, in terms of the spectrum of the underlying Liouville operator for deterministic systems or the Fokker-Planck operator for stochastic systems, but are in general difficult to estimate especially if the dimension of the state space is large. Section \ref{Sec_Theory_reduced_spec} below addresses the implication of estimating resonances from partial observations in a reduced state space, and how these reduced resonances relate to the full RP resonances.   

\subsection{Generalities}

We consider Stochastic Differential Equations (SDEs) of the form:
\begin{equation}\label{Eq_SDE}
\d X= F(X) \d t + D(X) \d W_t, \qquad X \in \mathbb{R}^d.
\end{equation}
Here $W_t=(W_t^1,\cdots,W_t^q)$ denotes an $\mathbb{R}^q$-valued Wiener process ($q$ not necessarily equal to $d$) whose components are mutually independent Brownian motions. 

In Eq.~\eqref{Eq_SDE}, the drift part  is provided by a (possibly nonlinear) vector field $F$ of $\mathbb{R}^d$, and the (also possibly nonlinear) stochastic diffusion in its It\^o version, given by $D(X) \d W_t $, has its $i^{th}$-component ($1\leq i\leq d$) given by

\begin{equation}
\bigg[D(X) \d W_t \bigg ]_i=\sum_{j=1}^q D_{ij} (X) \d W_t^j, \; \; q\geq 1.
\end{equation}
Throughout this article, we assume that the vector field $F$ and the (non-zero rank) matrix-valued function 
\begin{equation*}
D:\mathbb{R}^d \rightarrow \mbox{Mat}_{\mathbb{R}}(d\times q), 
\end{equation*}
satisfy regularity conditions that guarantee the existence and the uniqueness of mild solutions, as well as the continuity of the trajectories; see, e.g.~\cite{cerrai2001second,flandoli2010flow} for such conditions in the case of locally Lipschitz coefficients. Sometimes the resulting stochastic process 
solving Eq.~\eqref{Eq_SDE} emanating from $x$ at $t=0$, will be denoted by $X_t^x$, when dependence on the initial datum would be stressed, otherwise simply $X_t$.

It is well-known that the evolution of the probability density of the stochastic process, $X_t$, solving Eq.~\eqref{Eq_SDE},  is governed by the {\it Fokker-Planck equation}
\begin{equation}\label{Eq_FKP}
\partial_t  \rho(X,t) =\mathcal{A} \rho(X,t)= -\mbox{div}(\rho(X,t)F(X)) +\frac{1}{2} \mbox{div} \nabla (\mathbf{\Sigma}(X)  \rho(X,t)), \;\; X\in \mathbb{R}^d,
\end{equation}
with $\mathbf{\Sigma}(X)=D(X) D(X)^T$ denoting the diffusion tensor. In practice, one is interested in stationary solutions $\mu$ to Eq.~\eqref{Eq_FKP} that provide a statistical description of the asymptotic behavior of the solutions to Eq.~\eqref{Eq_SDE}, and thus of probability density functions (PDFs), typically  obtained  as marginal distributions of $\mu$.

What is less-known however is that the spectral properties of the 2nd-order differential operator, $\mathcal{A}$, inform about  fundamental objects such as the power spectra or correlation functions computed typically along a stochastic path of  Eq.~\eqref{Eq_SDE}.   To understand these relationships,  we derive in the next subsection, decomposition formulas of correlations and power spectra in terms of the spectral elements of an operator closely related to $\mathcal{A}$, namely the generator $K$ of the associated Markov semigroup.

\subsection{Ruelle-Pollicott (RP) resonances and the spectral decomposition of correlation functions}\label{Sec_spec_Pt}

As recalled in Appendix \ref{Sec_ergodic+mixing} (Theorem  \ref{Thm_cond_mixing}), the existence of an invariant measure $\mu$ ensures that the Markov semigroup $P_t$ (see Appendix \ref{Sec_Marov_smg}) associated with Eq.~\eqref{Eq_SDE}, 
is a strongly continuous semigroup in $L^p_\mu(\mathbb{R}^d)$, for any $p\geq 1$. This property alone allows us to get useful decomposition formulas of correlation functions 
$C_{f,g}(t)$ given by
\be\label{Cor_func_def}
C_{f,g}(t)=\int f\cdot P_t g \d \mu  - \int f \d \mu \int g \d \mu.
\ee
This is the content of Corollary \ref{Cor_decomp} below.  It makes use of elements of the spectral theory of strongly continuous semigroups that we recall below and apply to Markov semigroups.

In that respect,  recall that the {\it essential growth bound} of a semigroup $\mathcal{T}=(T(t))_{t\geq 0}$ is defined as 
\be\label{Def_ess-bound}
\omega_{ess}(\mathcal{T})=\underset{t>0}\inf \;  \frac{1}{t} \log \| T(t)\|_{ess}, 
\ee
where 
\be\label{Def_essnorm}
\| T(t)\|_{ess} =\inf \Big\{ \|T(t)- \mathcal{C}\|_{\mathcal{L}(\mathcal{E})} \, : \, \mathcal{C} \mbox{ is  a linear and compact operator of } \mathcal{E}=L^p_{\mu}(\mathbb{R}^d)\Big\}.
\ee
In other words, $\| T(t)\|_{ess}$ measures the distance of $T(t)$ to the set of linear and compact operators of $L^p_{\mu}(\mathbb{R}^d)$. 
In the case this distance approaches zero as $t \rightarrow \infty$, the strongly continuous semigroup  $\mathcal{T}$ is called {\it quasi-compact.}
The theory of quasi-compact strongly continuous semigroups shows that the latter property is equivalent to $-\infty<\omega_{ess}(\mathcal{T})<0$; see \cite[Prop.~V.3.5]{Engel_Nagel}. 

A semigroup is {\it eventually compact} if there exists $t_0>0$ such that $T(t_0)$ is compact. A semigroup that is eventually compact satisfies thus $\omega_{ess}(\mathcal{T})=-\infty$.  At the same time, the definition of \eqref{Def_ess-bound} allows for semigroups that are not eventually compact while their essential growth bound is negative infinity\footnote{For instance any semigroup $\mathcal{T}$ such that $\| T(t)\|_{ess}\leq M \exp{(-\epsilon t^{\alpha})}$,  with $\epsilon >0$ and $0<\alpha, M <1$.}.

In all cases, we have the following decomposition theorem of strongly continuous semigroups \cite[Theorem V.3.1]{Engel_Nagel} that we apply to the case of Markov semigroups possessing an invariant measure.
\begin{thm}\label{Thm_decomp}
Let $\mu$ be an  invariant measure of  a Markov semigroup $\mathcal{P}=(P_t)_{t\geq 0}$ and let $K$ be the corresponding generator in $L^2_\mu(\mathbb{R}^d)$. 
Let $\lambda_1,\cdots,\lambda_N$ be in $\sigma(K)$ satisfying $\Re \lambda_1,\cdots,\Re \lambda_N > \omega_{ess}(\mathcal{P})$.  Then the following properties hold:
\bi

\item[(i)] The $\lambda_j$ are isolated spectral values of $K$ with finite algebraic
multiplicity. 

\item[(ii)] If $\Pi_1,\cdots,\Pi_N$ denote the corresponding spectral projections
and $m_1,\cdots,m_N$ the corresponding orders of poles of the resolvent of $K$, then
\be\label{P-decomp1}
P_t =\sum_{j=1}^N T_j(t)  +R_N(t)
\ee
where 
\be\label{P-decomp2}
T_j(t)=\Big[\sum_{k=0}^{m_j-1} \frac{t^k}{k!} (K-\lambda_j \textrm{Id})^{k}\Big]  e^{\lambda_j t}\Pi_j,
\ee
and where for every $\epsilon >0$, there exists $M >0$ such that
\be\label{P-decomp3}
\| R_N(t)\|_{\mathcal{L}(L^2_\mu(\mathbb{R}^d))}  \leq M e^{(\omega^\ast_N+\epsilon) t}, \; \forall\, t\geq0,
\ee
with
\be \label{Est_omegaN}
\omega^\ast_N =\sup \{ \omega_{ess}(\mathcal{P})\} \cup \{\Re \lambda \; : \;  \lambda \in \sigma(K)\backslash \{\lambda_1,\cdots,\lambda_N\}\}.
\ee
\ei 
\end{thm} 
In what follows we denote by $\langle \cdot ,\cdot \rangle_{\mu}$ the inner product in $L^2_\mu$.
We have then the following corollary regarding the decomposition of correlation functions. 
 \bc\label{Cor_decomp}
 Let $\mu$ be an invariant measure of a Markov semigroup $\mathcal{P}=(P_t)_{t\geq 0}$ and let $K$ be the corresponding generator in $L^2_\mu$. 
 Let $f$ and $g$  in  $L^2_{\mu}$ such that $\langle f\rangle_{\mu}=\langle g\rangle_{\mu}=0$, then given the $\lambda_js$ and the $\Pi_js$ as defined in Theorem \ref{Thm_decomp}, the correlation function $C_{f,g}(t)$ associated with the observables  $f$ and $g$ has the following expansion 
\be\label{Eq_decomp_corr1}
\boxed{C_{f,g}(t)= \sum_{j=1}^N \sum_{k=0}^{m_j-1} \frac{t^k}{k!} e^{\lambda_j t} \bigg(\int f (x) (K-\lambda_j \textrm{Id})^{k}   (\Pi_j g) (x) \d \mu (x)\bigg) + \mathcal{Q}_N(t),}
\ee
with
\be
|\mathcal{Q}_N(t)|\leq M \|f\|_{L_{\mu}^2}\|g\|_{L_{\mu}^2} e^{(\omega_N^\ast +\epsilon) t}, \; \forall\, t\geq0. 
\ee
In the particular case of a discrete spectrum, $\sigma(K) =\{\lambda_1,\lambda_2,\cdots\}$, with $\Re \lambda_{j+1} \leq \Re \lambda_{j} $ for all $j\in\mathbb{N}^\ast$, and $\Re \lambda_j \rightarrow -\infty$,  if $m_j=1$ for all $j\geq 1$, then the expansion \eqref{Eq_decomp_corr1} takes the form:
\be\label{Eq_decomp_corr2}
\boxed{C_{f, g}(t)= \sum_{j=1}^{\infty} e^{\lambda_j t} \langle f, \psi_j \rangle_\mu \langle\psi_j^\ast, g \rangle_\mu,}
\ee
 where $\psi_j$ denotes the $L^2_{\mu}$-eigenfunction associated with $\lambda_j$ and  $\psi_j^\ast$, the associated $L^2_{\mu}$-eigenfunction of the adjoint operator $K^\ast$ of $K$.
\ec
 
 \begin{proof}
The decomposition \eqref{Eq_decomp_corr1} is a direct consequence of the definition of $C_{f,g}$  in \eqref{Cor_func_def} and of \eqref{P-decomp1}-\eqref{P-decomp3}. 
The proof of  \eqref{Eq_decomp_corr2} consists of noting that when $m_j=1$ in \eqref{Eq_decomp_corr1}, the integrals therein reduce to
\bea\label{Eq_coeff}
\int f (x) (\Pi_j g) (x) \d \mu (x)&=  \int f(x) \langle\psi_j^\ast, g \rangle_\mu \psi_j(x) \d \mu(x),\\
                                              &=\langle f, \psi_j \rangle_\mu \langle\psi_j^\ast, g \rangle_\mu,
\eea
and that $\mathcal{Q}_N(t)\underset{N\rightarrow \infty}\longrightarrow 0$ in virtue of our assumptions about $\sigma(K)$ and the reminder estimate \eqref{P-decomp3}. \qed
\end{proof}

\begin{figure}
\centering
	\begin{subfigure}{0.48\textwidth}
\centering	
\includegraphics[height=.8\textwidth, width=.85\textwidth]{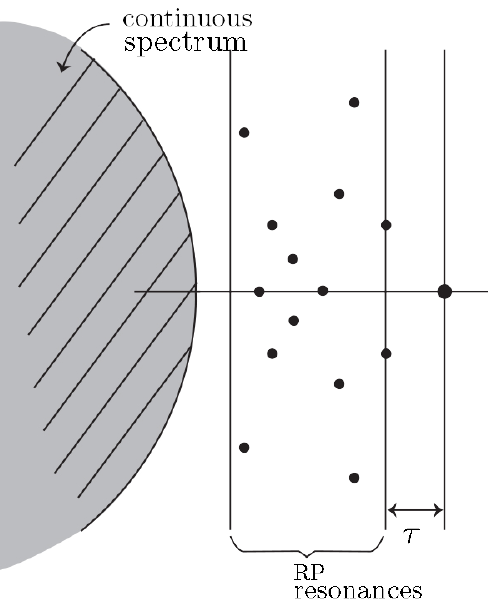}
\caption{Spectrum of Markov semigroup generator, $K$.}
	\end{subfigure}
		\begin{subfigure}{0.48\textwidth}
		\centering
		\includegraphics[width=.8\textwidth,height=.85\textwidth]{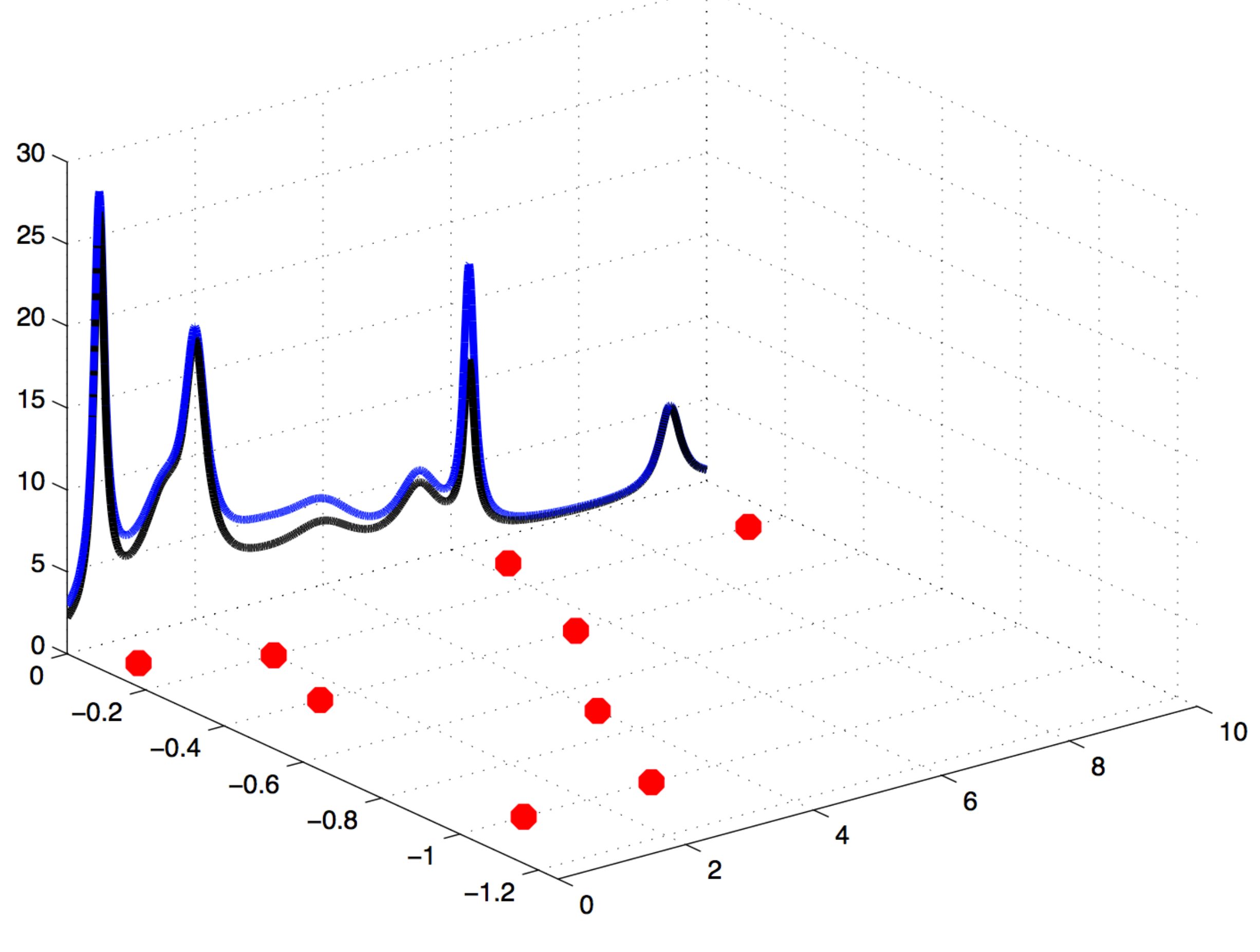}
		\caption{Correspondence between the PSD and RP resonances according to~\eqref{Eq_PSD_RPs}.}
	\end{subfigure}
	\caption{{\footnotesize{\bf Panel (a):} The Ruelle-Pollicott (RP)  resonances are shown here in the left half plane. 
The isolated eigenvalues of the Markov semigroup generator, $K$; they are represented by red dots in Panel (b) and by black dots here. The rightmost vertical line represents the imaginary axis above which the power spectrum  lies; see Panel (a) for another perspective.  The rate of decay of correlations is controlled by the spectral gap $\tau$; see Appendix \ref{Sec_spec_gap}. {\bf Panel (b):} [Courtesy of Maciej Zworski] The imaginary part of the RP resonances corresponds to the location of a peak in the PSD  (black curve lying above the imaginary axis) and the real part to its width. In blue is represented a reconstruction of the PSD based on RPs; a discrepancy is shown here to emphasize that in practice the RPs are very often only estimated/approximated; see~\cite{Chek_al14_RP}. }}\label{RP_schema}
\end{figure}

The $\lambda_j$ such that $\Re \lambda_j > \omega_{ess}(\mathcal{P})$, are called the {\bf RP resonances}. In other words they correspond to the point spectrum of $K$. Note that there exist stochastic processes for which $R_N=0$ (i.e.~no essential spectrum), for instance  Ornstein-Uhlenbeck processes \cite{metafune2002spectrum}. Remarks \ref{Rem_discrete_spec}-(ii)/(iii) point out other (nonlinear) stochastic processes that do not have an essential spectrum. Panel (a) of Fig.~\ref{RP_schema} shows a schematic of a general case, where the RP resonances are contained within a strip of the complex plane, away from the imaginary axis and the essential spectrum of $K$.  
The rate of decay of correlations is controlled by the spectral gap, $\tau=\mbox{gap}(K)$ defined in \eqref{Formula_gap}. Appendix \ref{Sec_spec_gap} reports on a broad class of  Markov semigroups (and thus SDEs) whose generator possesses a spectral gap; see also e.g.~\cite{froyland1997computer,Baladi2000,giulietti2013anosov,Chek_al14_RP} and references therein for results in the deterministic context.

 \br\label{Rmk_Dynkin}
\hspace*{2em}  \vspace*{-0.4em}
\bi

\item[(i)] {\bf Decay of correlations}. Under the conditions of Corollary \ref{Cor_decomp}, one obtains that $C_{f, g}(t) \underset{t\rightarrow \infty}\longrightarrow 0$, for any $f,g$ in $L^2_{\mu}$ such that $\int f \d \mu  =\int g\d \mu =0$, and without assuming $P_t$ to be strong Feller and irreducible.

\item[(ii)]  {\bf Transfer operator}. By working in the Hilbert space $L^2_{\mu}$ allows us to define the transfer operator $\mathfrak{L}_t $ acting on functions in $L^2_{\mu}$, by the duality formula:
\be \label{Eq_exchange2}
\int  f\cdot P_t g \d \mu  =  \int \mathfrak{L}_t f   \cdot  g \d \mu, \; f, g\in L^2_{\mu}, \; t\geq 0.
\ee
In a similar fashion that $P_t$ is related to the Kolmogorov operator $\mathcal{K}$  via its generator (see (iii) below),  the 
transfer operator  $\mathfrak{L}_t$ is related to the Fokker-Planck operator $\mathcal{A}$ defined in \eqref{Eq_FKP}. 

\item[(iii)]  {\bf RP resonances and Kolmogorov operator}. There are relationships between the abstract operator $K$ and the concrete Kolmogorov operator $\mathcal{K}$ defined in \eqref{Def_K0-additive}.   In that respect, a theorem of Dynkin ensures (e.g.~\cite[p.~258]{rogers2000diffusions}) that if $(P_t)_{t\geq0}$ is a {\it Feller semigroup} of a Markov process with continuous paths such that $C_0^{\infty}(\mathbb{R}^d)$
is contained in the domain $D(K)$ of the generator $K$, then there exist functions $a_{ij}$, $b_i$, and $c$  in $C(\mathbb{R}^d)$ $(1\leq i,j \leq d)$ such that for any $x$, the matrix $\big(a_{ij}(x)\big)_{1\leq i,j\leq d}$ is non-negative definite, $c(x) \leq 0$,%
\footnote{Furthermore if the process is non-explosive then $c\equiv 0.$ This excludes the cases for which the underlying Markov process leaving at time 0 from $x$ in $\mathbb{R}^d$ escapes to infinity at some finite time $t > 0$. This article is not concerned with explosive stochastic processes.}
 and
\be\label{Eq_K}
Kf (x) =\sum_{i,j=1}^d a_{ij}(x) \partial_{ij}^2 f + \sum_{i=1}^d b_i (x)\partial_i f + c(x) f, \; \; f\in C_0^{\infty}(\mathbb{R}^d).
\ee
Furthermore,
\be\label{Expression_aij}
b_i(x)=F_i(x),  \; \;\; a_{ij}(x)=\frac{1}{2}\big(D(x) D(x)^T\big)_{ij}.
\ee

For such reasons we sometimes refer to the RP resonances as the RP resonances of the Kolmogorov operator $\mathcal{K}$, or the  Fokker-Planck operator $\mathcal{A}$ (by duality). This language is often adopted in Part III, more oriented towards applications than the present Part I. 

We mention that the spectral analysis of general 2nd-order operators, typically non-self adjoint such as  given in \eqref{Eq_K}, is not an easy task in practice, especially when $p$ is large; see \cite{davies1996spectral,davies2007linear}.  Instead, {\mkr the notion of reduced RP resonances based on Corollary \ref{Cor_approx_RP} below} (see also Sec.~\ref{Sec_Theory_reduced_spec}) informs us about the spectral elements of {\mkr some coarse-grained version of  $\mathcal{K}$; see Theorem \ref{Thm_assymp_gene} and Remark \ref{Rmk_approx_RPs}-(iii).}
We refer nevertheless to \cite{eckmann2003spectral,herau2004isotropic} for a detailed (and enlightening) study of the ``shape'' of the spectrum (in $L^2_\mu$ for \cite{herau2004isotropic} and Sobolev spaces for \cite{eckmann2003spectral}) of a broad class of hypoelliptic operators such as arising in the theory of Langevin dynamics; see also \cite{ottobre2012exponential}.  
\ei
\er

 \br\label{Rem_discrete_spec}
\hspace*{2em}  \vspace*{-0.4em}
\bi
\item[(i)] Note that the sum in \eqref{Eq_decomp_corr2}  starts actually at $j=2$ since $\psi_1=\mathbbm{1}_{\mathbb{R}^d}$ (since $1\in \sigma(P_t)$) and $\langle f\rangle_{\mu}=0$, by assumption. 

\item[(ii)] When $\sigma(K) =\{\lambda_1,\lambda_2,\cdots\}$ with $\Re \lambda_{j+1} \leq \Re \lambda_{j} $ for all $j\in\mathbb{N}^\ast$, and $\Re \lambda_j \rightarrow -\infty$, is only assumed (without requiring $m_j=1$) then  \eqref{Eq_decomp_corr1} holds with $N=\infty$ and $\mathcal{Q}_N=0$.  Such an hypothesis on the spectrum of $K$ is systematically satisfied for the Markov semigroups that are eventually compact and possess an invariant measure $\mu$; see  \cite[Cor.~V.3.2]{Engel_Nagel}. A large subclass of such Markov semigroups are the compact ones. 
It includes a broad class of SDEs such as the following gradient dynamics on $\mathbb{R}^d$ subject to noise:
\be\label{SDE_gradient}
\d x=-\nabla V(x) \d t + \sqrt{2 \beta^{-1}}\d W_t, \; \beta>0,
\ee
where $W_t$ is a $d$-dimensional Wiener process. The diffusion is in that case elliptic, so it is straightforward that the assumptions of H\"ormander's theorem are satisfied, and the corresponding transition probabilities are smooth, provided that $V$ is smooth.  In fact for a broad class of smooth potentials $V$ that satisfy  appropriate growth conditions, the compactness of $P_t$ is ensured; e.g.~\cite[Theorem 8.5.3]{lorenzi2006analytical} and \cite[Lem.~1.2]{ganidis1999convergence}.  See also \cite{metafune2002compact} for conditions ensuring the compactness of Markov semigroups that are not necessarily associated with gradient dynamics. 

\item[(iii)]  Note that if $P_t$ is eventually compact and possesses a unique ergodic invariant measure $\mu$, then there exists $\alpha >0$ and $M\geq 1$ and 
\be\label{Eq_expo_conv2}
\| P_t \varphi -\langle \varphi\rangle_{\mu}\|_{L^2_\mu}\leq M e^{-\alpha t} \| \varphi -\langle \varphi\rangle_{\mu}\|_{L^2_\mu}.
\ee
Note that if $P_t$ has a regularizing effect sending $L_{\mu}^2$ into the Sobolev space $W_{\mu}^{1,2}$ in finite time and $W_{\mu}^{1,2}$ is compactly embedded into $L_{\mu}^2$ then the  Poincar\'e inequality holds \cite[Prop.~8.6.1]{lorenzi2006analytical} and $P_t$ is eventually compact.
The proof of \eqref{Eq_expo_conv2} is then a consequence of the theory of eventually compact semigroups \cite[Cor.~V.3.3]{Engel_Nagel} and the existence of a unique ergodic measure which ensures that $\mathbbm{1}_{\mathbb{R}^d}$ is the unique $L^2_{\mu}$-eigenfunction of $P_t$ (up to a multiplication factor) associated with the eigenvalue $1$.

\ei
\er

\subsection{Decomposition of the power spectrum}\label{Sec_PSD}
Given an observable $h:\mathbb{R}^d\rightarrow \mathbb{R}$ for the system $\eqref{Eq_SDE}$, we recall that  the {\it correlation spectrum} $S_h(f)$ is obtained by taking the Fourier transform of the correlation function $C_h (t)$, namely 

\begin{equation}\label{Eq_corr}
S_h(f) =\widehat{C_h}(f),  
\end{equation}
where $C_h$ given by \eqref{Cor_func_def} (for a given invariant measure $\mu$) with $f=g=h$, therein.

For a broad class of SDEs that possess an ergodic probability distribution $\mu$,   
the spectrum in $L^2_\mu$, $\sigma(K)$, of the Markov semigroup generator,  is typically contained in the left-half complex plane, $\{z\in \mathbb{C}\;:\: \Re(z)\leq 0\}$ and its resolvent  $R(z)=(z \mbox{Id}-K)^{-1}$, is a well-defined linear operator that satisfies

\begin{equation}\label{Eq_PSD_RPs}
S_h(f) =\int_{\mathbb{R}^d} h(X) \big[R(i f) h \big](X)\d \mu.
\end{equation}

In~\eqref{Eq_PSD_RPs}, the frequency $f$ lies in the complex plane $\mathbb{C}$, and the poles of the resolvent  $R(i f)$ \textemdash which correspond to the RP resonances \textemdash introduce singularities into $S_h(f)$. Once the power spectral density (PSD) is calculated, i.e.~once $|S_h(f)|$ is computed with $f$ taken to be real, these poles manifest themselves as peaks that stand out over a continuous background at the frequency $f$ if the corresponding RP resonances with imaginary part $f$ (or nearby) are close enough to the imaginary axis. The continuous background may have different origins.  In the case of a pure point spectrum, it is due to RP resonances located far from the imaginary axis.  In the presence of a continuous spectrum and RP resonances (such as shown in Panel (a) of Fig.~\ref{RP_schema}), if the latter are close to the imaginary axis, then the continuous background of the PSD is mainly due to the contribution of the continuous part of $\sigma(K)$ lying typically in a sector $\{z\in \mathbb{C}\;:\: \Re(z)\leq -\gamma\}$, for some $\gamma >0$.\footnote{While we recall that in such a case, the RP resonances are the  isolated eigenvalues of finite multiplicity,  lying within a strip $-\gamma < \Re(z) \leq 0$; see Panel (a) of Fig.~\ref{RP_schema}.} 

Formula \eqref{Eq_PSD_RPs} establishes furthermore a useful correspondence between PSD and RP resonances. 
Indeed, from  \eqref{Eq_PSD_RPs}   we infer that the imaginary part of the RP resonances corresponds to the location of a peak in the PSD  and the  absolute value of the real part to its width; see Panel (b) of Fig.~\ref{RP_schema}. 

Let us take $f=g=h$ in \eqref{Eq_decomp_corr1}. By denoting by $\alpha_j^k(h)$ the coefficients resulting as integration with respect  to $
 \mu$ in  \eqref{Eq_decomp_corr1}, the latter decomposition formula
writes
\begin{equation}\label{Eq_decomp_corr1_b}
C_{h}(t)= \sum_{j=1}^{N} \bigg(\sum_{k=0}^{m_j-1} \frac{t^k}{k!}  \alpha_j^k(h) \bigg)e^{\lambda_j t} + \mathcal{Q}_N(t),
\end{equation}  
where $\mathcal{Q}_N(t)$ exhibits typically a decay property associated with properties of the essential spectrum of $\mathcal{A}$.  Note that  the $\lambda_j$ do not depend on the observable $h$, but that the $\alpha_j^k(h)$ do. 

If we assume that $\Re(\lambda_j) < 0$ for $j > 0$, that each $\lambda_j$ is simple ($m_j=1$) and the absence of an essential spectrum for $K$,  then 
the correlation $C_{h}(t)$ in \eqref{Eq_decomp_corr1} takes the simpler form of a weighted sum of complex exponentials (i.e.~\eqref{Eq_decomp_corr2}), where $\alpha_j^k=\alpha_j$ as given by \eqref{Eq_coeff},  and the corresponding correlation spectrum $S_{h}(f)$ possesses itself a similar decomposition in terms of Lorentzian functions, namely:
\be
	S_{h}(f) = - \frac{1}{\pi} \sum_{j = 1}^\infty \alpha_j(h) \frac{\Re(\lambda_j)}{(f - \Im(\lambda_j))^2 + \Re(\lambda_j)^2}, \qquad f\in \mathbb{C}. 
	\label{eq:spectralPower}
\ee
It is noteworthy that such Lorentzian decompositions of the PSD are frequently encountered in applications such as for instance in spectroscopy due to the presence of resonances; e.g.~\cite{gruberbauer2009detection,kallinger2010asteroseismology,nonnenmacher2015decay}. 
Note that RP resonances appear also in the decomposition of the susceptibility function arising in linear response theory \cite[Eq.~(5)]{lucarini2018revising}; see  
\cite{ruelle2005differentiating,cessac2007does,ruelle2009review}. 

In summary, the decompositions~\eqref{Eq_decomp_corr1} and~\eqref{eq:spectralPower} inform us about the following features:
\begin{itemize}
	\item[(i)] Each RP resonance is associated with an exponential contribution to the decay of correlation.
	\item[(ii)] The closer an eigenvalue to imaginary axis, the slower the decay.
	\item[(iii)] In the limit of purely imaginary eigenvalues, the associated contributions
	to the correlation functions are purely oscillatory and prevent the decay of correlations.
	\item[(iv)] The angular frequency at which each contribution oscillates
	is given by the imaginary part of the associated eigenvalue.
	\item[(v)] Eigenvalues close to the imaginary axis are associated with resonances (i.e.~peaks) in the
	power spectrum.
	The spectral peak is located at the frequency given by the imaginary part of the eigenvalue
	and its width is proportional to the absolute value of the real part.
	\item[(vi)] The contribution of each eigenvalue to a correlation function or a power spectrum
	is weighted as in \eqref{Eq_coeff}, corresponding to the projection
	of the observables $h$ onto the eigenfunctions of $K$ and its adjoint.
\end{itemize}

\section{Reduced Ruelle-Pollicott resonances}\label{Sec_Theory_reduced_spec}
We present in this section the main results that serve us as a foundational basis for the applications discussed hereafter and in Part III \cite{PartIII}. 
The goal is to characterize the resonances that can be extracted from low-dimensional observations in a reduced state space and how these resonances relate
to the RP resonances introduced in Sec.~\ref{Sec_spec_Pt} above.

The results of this section rely substantially on the general disintegration theorem of probability measures\footnote{Variation of this theorem is used in the study of spectral gaps for deterministic maps and is known as Rokhlin's disintegration theorem; see \cite{galatolo2020spectral}.}; see e.g.~\cite[p.~78]{dellacherie1978probabilities}. See also \cite[Theorem 5.3.1]{ambrosio2008gradient} and \cite{chang1997conditioning}.
The disintegration theorem states that given a probability measure $\mu$ on $\mathbb{R}^d$,  a finite-dimensional Euclidean space $V$ for which $1\leq \,$dim$(\V)<d$,  and a map $h: \mathbb{R}^d \rightarrow  V$ (Borel-measurable), then there exists a uniquely determined family of probability measures $\{\mu_v\}_{v\in V}$ such that, for $\m$-almost  all \footnote{i.e.~up to an exceptional set of null measure with respect to $\m$.} $v$ in $V$, $\mu_v$ is concentrated on the pre-image $h^{-1} (\{v\})$ of $v$, i.e.~$\mu_v \left(\mathbb{R}^d\setminus h^{-1} (\{v\}) \right) = 0,$ 
and such that for every Borel-measurable function $\phi: \mathbb{R}^d \rightarrow \mathbb{R}$,
\be\label{Eq_desint}
\int_{\mathbb{R}^d} \phi (x) \d \mu (x) = \int_{V} \Big(\int_{x\in h^{-1} (\{v\})} \phi(x) \d\mu_{v} (x) \Big)\d \m (v).
\ee
Here $\m$ denotes the push-forward in $V$ of the measure $\mu$  in $\mathbb{R}^d$, by the map $h$.

Hereafter, we apply this result for the proof of Theorem \ref{Thm_assymp_gene} when the reduced state space, $V$, is a subspace of  $\mathbb{R}^d$, and the mapping $h$ is the projector $\Pi_V$ onto $V$. In this case, a decomposition analogous  to \eqref{Eq_desint} holds for the measure $\mu$ itself, namely  
\be\label{Eq_desint2}
\mu(B\times F)=\int_{F} \mu_{v} (F) \d \m (\xi), \qquad B\times F \in \mathcal{B}(V)\otimes\mathcal{B}(W),
\ee
with $W$ being the subspace such that $V\oplus W=\mathbb{R}^d$. For later notation convenience, we denote by $\cX_v$ the pre-image $h^{-1} (\{v\})$ (for any $h$ surjective), namely 
\be\label{Def_Xv}
\cX_v=h^{-1} (\{v\}), \;\; \forall \; v \in V.
\ee

\subsection{Reduced Markov operators from partial observations of stochastic systems}
We start with Theorem \ref{Theorem_A_stoch} below that extends  \cite[Theorem A]{Chek_al14_RP} within the stochastic context of this article.     {\it Mutatis mutandis},  Theorem \ref{Theorem_A_stoch} as \cite[Theorem A]{Chek_al14_RP}, shows that --- given a reduced state space $\V$, a continuous map $h$ from $\mathbb{R}^d$ to $\V$, and a Markov semigroup $(P_t)_{t\geq 0}$ possessing an invariant measure on $\mathbb{R}^d$ --- a family of reduced Markov operators (acting on functions of $\V$) can be naturally associated with $(P_t)_{t\geq 0}$, $\V$ and $h$. This family characterizes a coarse-graining in the reduced state space $\V$ --- such as induced by the map $h$ --- of the actual transition probabilities associated with $(P_t)_{t\geq 0}$ that take place within the full state space $\mathbb{R}^d$.   Theorem \ref{Theorem_A_stoch} can be viewed as providing the analogue, in the SDE context dealt with in this article, of coarse-grained Markov representations used  for describing coarse-grained dynamics of Hamiltonian systems \cite{schutte1999direct} or more general Markov state models \cite{schutte2013metastability}.
  We emphasize that, as for these contexts, empirical estimates of the transition probabilities appearing in \eqref{Eq_main} of Theorem \ref{Theorem_A_stoch} are important for applications, as briefly discussed in Section \ref{Sec_reducedRP} below and in more details in Part III \cite{PartIII}.

\bt\label{Theorem_A_stoch}
Let $(P_t)_{t\geq 0}$ be a Markov semigroup that possesses an invariant measure $\mu$, and let $V$ be a reduced state space  for which $1\leq \mbox{dim}(V) < p$. Let $h:\mathbb{R}^d \rightarrow V$ be a continuous surjective function and let us denote by $\m=h_{\ast} \mu$, the push-forward of the measure $\mu$ by $h$.

Then there exists a time-dependent family of Markov operators $T_t$  acting on $L^1_{\m}(V)$ such that, for any Borel sets $B$  and $C$ of $V$, and any $t>0$,
\be\label{Eq_main} 
\frac{\langle T_t  \mathbbm{1}_B, \mathbbm{1}_C\rangle_{1,\infty}}{\m(B)}=\textrm{Prob}(S(t, \cdot) x\in h^{-1}(C) | x \in h^{-1}(B)),
\ee  
where $S(t,\omega)$ denotes the stochastic flow associated with $P_t$, and  $\mathbbm{1}_E$ denotes the characteristic function of a Borel set $E$.

Furthermore, if $\mu$ is ergodic then for any $f$  in $L^2_{\m}(V)$, it holds
\be\label{weak_ergodic}
\frac{1}{t} \int_0^t \int_{V} T_s f (v) \d \m (v) \d s \underset{t\rightarrow \infty}\longrightarrow\int_{V} f(v) \d \m(v).
\ee 
\et

\bp[Proof.] 

Let us introduce for any $t>0$ the operator $T_t$ which maps $f$ in $L_{\m}^1(V)$ to $T_t f$ given by
\be \label{Def_Tt}
T_t f(v)=\int_{x\in \cX_v} P_t [f \circ h](x) \d \mu_v(x),
\ee
where $\cX_v$ denotes the pre-image (in $\mathbb{R}^d$) of $v$ by $h$, and $\mu_v$ denotes the disintegrated probability measure associated with $h$, such as given by the disintegration theorem recalled above. The bracket-notation, $[\cdot]$, is used here in \eqref{Def_Tt} to emphasize that the formula should be read as follows: $P_t$ is first applied on $f \circ h$, then the integration with respect to $\mu_v$, is undertaken. Note that $T_t$ acts on observables of $V$, i.e.~on functions $\phi(v)$, with $v$ lying within the reduced state space. 

By integrating \eqref{Def_Tt} with respect to the measure $\m$, and applying \eqref{Eq_desint} with $\phi=P_t [f\circ h]$, one notes that 
\be
\int_{V}T_t f(v) \d \m (v)=\int_{\mathbb{R}^d} P_t [f\circ h] (x) \d \mu(x).
\ee
 
This last identity allows us to get the following estimates 
\bea 
\Big|\int_{V}T_t f(v) \d \m (v)\Big|& \leq \int_{\mathbb{R}^d} | P_t [f \circ h](x) |\d \mu (x)\\
 &\leq \int_{\mathbb{R}^d} |f \circ h (x)| \d \mu (x),
 \eea
the last inequality resulting from \eqref{Eq_Holder} applied with $p=1$ and the invariance of $\mu$. We have thus the estimate 
\be 
 \|T_t f \|_{L^1_{\m}(V)} \leq \int_{V} |f(v)| \int_{\cX_v}  \mathbbm{1}_{V} \d \mu_v  \d \m (v)=  \|f\|_{L^1_{\m}(V)},
\ee
which shows that $T_t$ maps $L^1_{\m}(V)$ into itself, i.e.~$T_t$ is a Markov operator.

The rest of the proof is also straightforward, and consists of noting that for any Borel sets $B$ and $C$ in $V$, one has 
\bea
\langle T_t  \mathbbm{1}_B, \mathbbm{1}_C\rangle_{L^1_{\m},L^{\infty}_{\m}} =&\int_{V} \mathbbm{1}_C(v) \Big[\int_{\cX_v} P_t[\mathbbm{1}_B \circ h](y) \d \mu_v(y)\Big]\d \m (v)\\
&= \int_{\mathbb{R}^d} \mathbbm{1}_{h^{-1}(C)}(x)P_t[\mathbbm{1}_B \circ h](x) \d \mu(x)\\
& = \int_{\mathbb{R}^d} \mathfrak{L}_t \mathbbm{1}_{h^{-1}(C)} \cdot \mathbbm{1}_B \circ h \, \d \mu\\
& = \int_{\mathbb{R}^d} \mathfrak{L}_t \mathbbm{1}_{h^{-1}(C)} \cdot \mathbbm{1}_{h^{-1}(B)} \,  \d \mu
\eea

Since $\m=h_{\ast }\mu$,  we deduce that $\m(B)=\mu(h^{-1}(B))$, and thus by dividing by $\m(B)$, one obtains that
\be
\frac{\langle T_t  \mathbbm{1}_B, \mathbbm{1}_C\rangle_{1,\infty}}{\m(B)}= \mbox{Prob}(S(t, \cdot) x\in h^{-1}(C) | x \in h^{-1}(B)).
	\quad B, C \in \mathcal{B}(V).
\ee

Finally the proof of \eqref{weak_ergodic} consists of applying again \eqref{Eq_desint} the characterization (iii) of the ergodicity recalled in
Definition \ref{Def_ergodic}. \qed
\ep

\br \label{Rmk_spectral-issue}
For any $f$  and $g$ in $L^2_{\m}(V)$, one can define correlation functions (in the reduced state space $V$) associated with the family of Markov operators $T_t$ by
\be\label{Cor_func2}
\widetilde{C}_{f,g}(t)=\int_V f\cdot T_t g \d \m  - \int_V f \d \m \int_V g \d \m.
\ee
These correlations inherit a nice property resulting from the disintegration formula \eqref{Eq_desint}, namely 
\be\label{Cor_func}
\widetilde{C}_{f,g}(t)=C_{f\circ h,g\circ h}(t),
\ee
where $C_{f\circ h,g\circ h}(t)$ is given by \eqref{Cor_func_def} in which $f$ (resp.~$g$) is replaced by $f\circ h$ (resp.~$g\circ h$). In particular, the decomposition formulas of Corollary \ref{Cor_decomp} apply to 
the observables $f\circ h$ and $g\circ h$, and provide in turn a spectral decomposition of correlation functions in the reduced state space $V$, in terms of the spectral elements of the generator $K$ of the Markov semigroup $(P_t)_{t \geq 0}$.

The issue, however, is that such a decomposition involves a spectral analysis of $K$ that is  
typically out of reach numerically when the dimension of the state space is large, and is difficult analytically in the general case. One would like instead to have at our disposal a decomposition in terms of the spectral elements associated with the reduced Markov operators $(T_t)_{t\geq0}$ rather than with the full Markov semigroup $(P_t)_{t\geq0}$. The stumbling block is that $(T_t)_{t\geq0}$ does not form a semigroup in general i.e.~
\be\label{non-Markov}
T_{t+s}\neq T_t T_s,
\ee
(see \cite{tantet2015early}) and one cannot rely on the spectral theory of semigroups such as used  in Section~\ref{Sec_spec_Pt}, to reach this goal.  
The next sections analyzes what type of useful spectral information may still be extracted from $T_t$. 
\er

\subsection{Pseudo generators of the family of reduced Markov operators and asymptotic behavior}\label{Sec_cond_exp}
We begin with the following Lemma that characterizes the pseudo-generators $G_t$ that can be associated with the family of reduced Markov operators $(T_t)_{t\geq 0}$, although the latter family does not form a semigroup in general. For related material, we refer to the recent works \cite{bittracher2015pseudogenerators,bittracher2015pseudo} concerned with the metastability analysis of Langevin dynamics. In a certain sense, Theorem \ref{Thm_assymp_gene} along with the {\mkr notion of reduced RP resonances based on Corollary \ref{Cor_approx_RP} below (see also Sec.~\ref{Sec_reducedRP}),} provide a complementary approach for the analysis of reduced Markov operators for situations not necessarily limited to Langevin dynamics. Theorem \ref{Thm_assymp_gene} shows in particular that the asymptotic behavior of $G_t$ as $t\rightarrow \infty$ (in a weak sense) is governed by the generator of a reduced SDE in which the dependence on the unobserved variables (lying outside of the reduced state space $\V$) has been averaged out; see \eqref{Eq_G_Markov} below.  In other words, Theorem \ref{Thm_assymp_gene} provides a useful relationship between the reduced non-Markovian process associated with the family $(T_t)_{t\geq 0}$ and the Markov process associated with the conditional expectation induced by the observable $h$.

\begin{lemma}\label{Lem_pseudo}
Let $K$ denote the generator of $P_t$ in $L^2_\mu(\mathbb{R}^d)$ with domain $D(K)$ as defined in \eqref{Def_domA}. Let $h:\mathbb{R}^d \rightarrow V$ be a continuous surjective observable.
Then, the Markov operator $T_t$ of Theorem \ref{Theorem_A_stoch} possesses for each time $t$, a pseudo-generator, i.e.~for any $f$ in $L^2_{\m}(V)$ for which $f\circ h$ lies in  $D(K)$, $\underset{s\rightarrow 0}\lim (T_{t+s}f-T_tf)/s$ exists and is given by 
\be\label{Eq_Gt}
G_t f (v)= \int_{\cX_v} P_t K[f\circ h] (x) \d \mu_v(x),
\ee
whose domain is thus 
\be\label{dom_X}
D(G_t)=\{f \in L^2_{\m}(V) \, :\, f\circ h \in D(K)\}. 
\ee
\end{lemma}

\bp[Proof.]
The proof is elementary and consists of noting that for any $f$ in  $L^2_{\m}(V)$ for which $f\circ h$ lies in $D(K)$, we have
\be
 \frac{T_{t+s}f-T_tf}{s}=\int_{\cX_v}\frac{P_{t+s}-P_t}{s} [f\circ h] (x) \d \mu_v(x)  \underset{s\rightarrow 0}\longrightarrow \int_{\cX_v} P_t K[f\circ h] (x) \d \mu_v(x),
\ee 
from the basic properties of a strongly continuous semigroup such as applied to $P_t$ (e.g.~\cite[Lem.~II.1.3]{Engel_Nagel}) and the dominated convergence theorem. \qed
\ep

\begin{thm}\label{Thm_assymp_gene}
Assume that the assumptions of Theorem \ref{Theorem_A_stoch} hold. If $\mu$ is ergodic, then there  exists a closed linear operator $\mathcal{G}$ acting on $L^2_{\m} (V)$ such that for any $f$  in $L^2_{\m}(V)$
\be\label{Eq_assymp_Gt}
\frac{1}{t}\int_0^t\int_{V} G_t f (v) \d \m (v) \underset{t\rightarrow \infty}\longrightarrow \int_{V} \mathcal{G} f (v) \d \m (v), 
\ee
and for which $\mathcal{G}$ is dissipative in the sense that
\be\label{Eq_dissip_G}
\Re \langle \mathcal{G} f,f \rangle_{L^2_{\m}(V)} \leq 0, \; \forall \, f \in D(\mathcal{G}).
\ee

Furthermore, if $\mu$ is strongly mixing then
\be\label{Eq_assymp_Gt_strong-mixing}
\lim_{t \rightarrow \infty}  \int_{V} G_t f (v) \d \m (v)=\int_{V}\mathcal{G} f (v) \d \m (v).
\ee

Recall that the Markov semigroup $P_t$ from which $G_t$ is defined, is associated to Eq.~\eqref{Eq_SDE}.
If $h$ is a projector, the operator $\mathcal{G}$ is densely defined, and if $P_t$ is a Feller semigroup for which  $C_0^{\infty}(\mathbb{R}^d) \subset D(K)$, then $\mathcal{G}$ possesses the following differential expression 
\be\label{Eq_G}
\G f (v) =\frac{1}{2}\sum_{i,j=1}^{{\rm dim}(\V)} \overline{\mathbf{\Sigma}}_{ij}(v) \partial_{ij}^2 f + \sum_{i=1}^{{\rm dim}(\V)} \overline{F}_i (v)\partial_i f , \; \;v\in V, \; \;  f\in C_0^{2}(\V).
\ee
where 
\be\label{coeff_G}
\overline{\mathbf{\Sigma}}_{ij}(v)=\int_{\cX_v} \mathbf{\Sigma}_{ij}(x) \d \mu_v (x), \;\; \overline{F}_{i}(v)=\int_{\cX_v} F_{i}(x) \d \mu_v (x), \;\;  x\in\mathbb{R}^d, \; v \in \V,
\ee
where $F$ and $\mathbf{\Sigma}$ denote respectively the drift part and the diffusion tensor associated with Eq.~\eqref{Eq_SDE}.

In other words, $\G$ is the generator of the Markov process associated with the {\bf conditional expectation} induced by $h$,
namely with a reduced  (``effective'') SDE in $\V$ of the form:
\be\label{Eq_G_Markov}
\boxed{
\d v= \overline{F}(v) \d t+\sigma(v) \d W_t^{\V}, \;\; v \in \V,
}
\ee  
with $W_t^{\V}$ denoting a Brownian motion in $\V$ and for $1\leq i,j\leq \mbox{dim}(\V)$,
\be\label{Eq_diffusion_averaged}
\overline{\mathbf{\Sigma}}_{ij}(v)=(\sigma(v)\sigma(v)^{\textrm{T}})_{ij}= \sum_{k=1}^q \overline{D_{ik}(x)D_{jk}(x)}, \; \; x\in \mathbb{R}^d, \; v \in \V,
\ee
where the $D_{\ell k}$ are the diffusion coefficients of the original SDE \eqref{Eq_SDE}, and $\overline{(\cdot)}$ denotes the averaging over the disintegrated measure $\mu_v$.

\end{thm}

\bp[Proof.]

\noindent {\bf Step 1: Proof of \eqref{Eq_assymp_Gt} and \eqref{Eq_dissip_G}.}
Let $f$ be in $D(G_t)$ given in \eqref{dom_X}. Then 
\bea
\int_{V} G_t f (v) \d \m (v) &= \int_{V}\int_{\cX_v} P_t K[f\circ h] (x) \d \mu_v(x) \d \m(v),\\
&=\int_{\mathbb{R}^d} P_t K[f\circ h] (x) \d \mu(x),
\eea
by application of the disintegration formula \eqref{Eq_desint}. 

Assuming the measure $\mu$ to be ergodic, the characterization (iii) of 
Definition \ref{Def_ergodic} allows us to infer \eqref{Eq_assymp_Gt} with 
$\mathcal{G}$ given by
\be\label{Eq_defG}
\mathcal{G} f (v) =\int_{\cX_v} K[f\circ h] (x) \d \mu_v (x),
\ee
by application once more of the disintegration formula. 

The same formula ensures furthermore that
\be
\int_{V} \mathcal{G} f \cdot f \d \m = \int_{\mathbb{R}^d} K [f\circ h] \cdot [f\circ h] \d \mu
\ee
and by taking the real part,  \eqref{Eq_dissip_G} follows from \eqref{Eq_dissip_K}.

\vspace{1ex}
\noindent {\bf Step 2: $\mathcal{G}$ is closed.} Let $\{f_n\}$ be a sequence in  $D(\mathcal{G})$ converging to $f$ in $L_{\m}^2(V)$
such that $\G f_n \rightarrow g $ in $L_{\m}^2(V)$ as $n\rightarrow \infty$. 

Then from \eqref{Eq_defG}, one has
\bea\label{Eq_step_closed1} 
\int_{\V}  |\G f_n -g|^2 \d \m= -2\int_{\mathbb{R}^d} & K [f_n\circ h] (x) \cdot (g \circ h (x)) \d \mu + \\
&\int_{\mathbb{R}^d} |K [f_n\circ h] (x) |^2\d \mu  + \int_{\mathbb{R}^d} |g\circ h (x)|^2 \d \mu(x).\\
\eea
Now since $K$ is closed in $L^2_{\mu}(\mathbb{R}^d)$,  $K [f_n\circ h] \underset{n\rightarrow\infty}\longrightarrow K [f\circ h]$ and $f \circ h \in D(K),$ namely
$f \in D(\G)$. This shows furthermore that one can pass to the limit in \eqref{Eq_step_closed1}, which gives
\be\label{Eq_step_closed2}
2\int_{\mathbb{R}^d} K [f\circ h] (x) \cdot (g \circ h (x)) \d \mu=\int_{\mathbb{R}^d} |K [f\circ h] (x) |^2\d \mu+\int_{\mathbb{R}^d} |g\circ h (x)|^2 \d \mu(x),
\ee
 The identity \eqref{Eq_step_closed2} says nothing else than 
\be
\|K [f\circ h] -g\circ h \|^2_{L^2_{\mu}}=0,
\ee
i.e.~$K [f\circ h] =g\circ h$, $\mu$-almost everywhere, which gives
\be
\int_{\cX_v} K [f\circ h] (x) \d \mu_v(x) =\int_{\cX_v} g\circ h (x) \d \mu_v(x)=g(v), \; \mbox{ for  } \m\mbox{-almost every } v \mbox{ in } \V. 
\ee
We have thus proved $\G f =g$, with $f$ in $D(K)$, i.e.~that $\G$  is closed.

\vspace{1ex}
\noindent {\bf Step 3: $\mathcal{G}$ is densely defined if $h$ is a projector.}  The purpose is to prove that $D(\G)$ is dense in $L^2_{\m}(\V)$. 
Assume by contradiction that there exist $g$ in $L^2_{\m}(\V)$ for which any sequence in $D(\G)$ would stay at a certain distance (in $L^2_{\m}(\V)$) from $g$ or would have $g$ as a point of accumulation, but not unique. 

{\mkr Let us take $V=\mathbb{R}^q$ with $1 \leq q<d$ (since $h$ is a projector) and define $G$ the mapping from $V$ to $V$ such that $\Pi_j G=g$, for all $1\leq j\leq q$. Since $D(K)$ is dense in $L^2_{\mu}(\mathbb{R}^d)$ (\cite[Cor.~2.5 p.~5]{Pazy1983}) there exists a sequence $\{\varphi_n\}$  in $(D(K))^d$ that converges to $G\circ h$ in $(L^2_{\mu}(\mathbb{R}^d))^d$,  such that in particular $\|\varphi_n-G\circ h\|$ (norm in $\mathbb{R}^d$) converges to zero in $L^2_\mu$.  
Let us define, almost everywhere (a.e.), a sequence $\{f_n\}$ in $(D(\mathcal{G}))^q$  from the semiconjugacy relation\footnote{Note that $f_n$ defined by \eqref{Eq_semi} implies that $\Pi_j f_n\circ h=\Pi_j \varphi_n$ belongs to $D(K)$ by construction, and thus $\Pi_j f_n$ belongs to $D(\mathcal{G})$ for every $1\leq j\leq q$, since the RHS of \eqref{dom_X} is also the domain of $\mathcal{G}$ as \eqref{dom_X} is independent on $t$.} 
\be\label{Eq_semi}
f_n\circ h=h\circ \varphi_n,
\ee
i.e.~such that the following diagram is commutative a.e.
\bes
\begin{CD}
  \mathbb{R}^d   @>\varphi_n \in (D(K))^d>>  \mathbb{R}^d\\
@VVhV        @VVhV\\
\V     @>f_n \in (D(\mathcal{G}))^q>>  \V
\end{CD}
\ees
Furthermore, since $h$ is a projector
\be
\int \|h(\varphi_n(x))-h(G\circ h (x)) \|^2 \d \mu (x)\leq \int \|\varphi_n(x)-G\circ h (x) \|^2 \d \mu (x)
\ee
which, thanks to \eqref{Eq_semi}, is equivalent to
\be\label{eq_interm}
\int \|f_n\circ h(x)-G\circ h (x) \|^2 \d \mu (x) \leq \int \|\varphi_n(x)-G\circ h (x) \|^2 \d \mu (x).
\ee
On the other hand, by application of the disintegration formula
\be
\int \bigg|\Pi_j \bigg(f_n\circ h(x)-G\circ h (x)\bigg) \bigg|^2 \d \mu (x) =  \int_{\V} |\Pi_j f_n(v) -g(v)|^2 \d \m (v), \; \; 1\leq j \leq q.
\ee
We have thus proved, thanks to \eqref{eq_interm} that for any $1\leq j \leq q$,
\be
\int_{\V} |\Pi_j f_n(v) -g(v)|^2 \d \m (v) \leq \int \|\varphi_n(x)-G\circ h (x) \|^2 \d \mu (x).
\ee
Now since $\|\varphi_n-G\circ h\|$ converges to zero in $L^2_\mu$ by construction, we conclude that $\Pi_j f_n$, which belongs to $D(\mathcal{G}$) (also by construction),  converges to $g$ in $L^2_{\m}(V)$, leading thus to a contradiction. 
}

\vspace{1ex}
\noindent {\bf Step 4: $\mathcal{G}$ provides the conditional expectation.} 
This step is a straightforward consequence of the Dynkin theorem (see Remark \ref{Rmk_Dynkin}-(iii)) and the representation formula \eqref{Eq_defG} which leads to the expression \eqref{Eq_G} of $\mathcal{G}$, by integration with respect to $\mu_v$. \qed

\ep

{\mkr We have then the immediate corollary.
\bc\label{Cor_approx_RP}
Under the assumptions of Theorem \ref{Thm_assymp_gene},  if $\mu$ is strongly mixing, the eigenvalues of $G_t$ provide approximations of the weak eigenvalues of  $\mathcal{G}$ in the sense that if $(\lambda_k^t,\psi_k^t)$ denotes an eigenpair of $G_t$, then 
\be\label{Eq_eigen_weak}
\int_V \mathcal{G} \psi^t_k \d \m  = \lambda^t_k\int_V \psi_k^t \d \m +\epsilon(t),
\ee
with $\epsilon(t) \rightarrow 0$, as $t \rightarrow \infty$.

Due to the definition of $G_t$ in \eqref{Eq_Gt} and of $\mathcal{G}$ in \eqref{Eq_defG}, we have also that 
\eqref{Eq_eigen_weak} holds with $\epsilon(t) \rightarrow 0$, as $t\rightarrow 0$.
\ec

The eigenvalues of $G_t$ in Corollary \ref{Cor_approx_RP} are called  that the {\bf reduced RP resonances} for the time-lag $t$. 
}

 \br\label{Rmk_approx_RPs}
\hspace*{2em}  \vspace*{-0.4em}
\bi

\item[(i)] {\bf Approximation of RP resonances}.
{\mkr In the context of homogenization of SDEs, explicit formulas for the generator $\mathcal{G}$ of the conditional expectation are available; see 
e.g.~\cite[Chaps.~11 \& 18]{pavliotis2008multiscale}. Within this context and adopting the language of this article, the estimates \cite[Eqns.~(5.12)-(5.13)]{crommelin2011diffusion} show that 
 the RP resonances (resp.~the eigenfunctions in $L^2_\mu$) of the full generator $K$ are approximated by the eigenvalues (resp.~the eigenfunctions in $L^2_\m$ after multiplication by the $\mu_x$-density) of $\mathcal{G}$, when the time-scale separation between the observed and unobserved variables is sufficiently large.    Thus, for such situations,  the reduced RP resonances provide (weak) approximations to the  actual RP resonances, due to Corollary \ref{Cor_approx_RP}. 

In the general case, error estimates between the RP resonances and the reduced RP resonances are difficult to derive. We refer though to Sec.~\ref{Sec_Langevin} for an example for which  the reduced RP resonances allow for very good approximation of correlation functions, even when the time-scale separation is weak.  There, it is the decomposition \eqref{Eq_decomp_corr2} that is approximated when the RP resonances (and corresponding modes) therein are replaced by the spectral elements of $G_t$.}

\item[(ii)] {\bf $\mathcal{G}$ describes long-  and short-time behaviors in $V$.} Due to \eqref{non-Markov}, the stochastic process associated with $T_t$ in the reduced state space $V$ is in general non-Markovian. Theorem \ref{Thm_assymp_gene} shows however that $\mathcal{G}$ plays an important role in the description of the asymptotic behavior of this stochastic process. 
{\mkr As pointed out in Corollary \ref{Cor_approx_RP}, the operator} $\mathcal{G}$ as defined in \eqref{Eq_defG} is exactly the pseudo-generator $G_t$ of $T_t$, as $t\rightarrow 0$. 

The operator  $\mathcal{G}$ relates thus the long time assymptotics ($t\rightarrow \infty$) of the reduced stochastic process to its infinitesimal characteristics which describe the short-time behavior ($t\rightarrow 0$), as for time homogeneous Markov processes.
We refer to \cite{legoll2010effective,legoll2017pathwise,zhang2017effective,duong2018quantification} for error estimates between the coarse-grained dynamics and the effective dynamics involving conditional expectation, in the context of Langevin or overdamped Langevin dynamics.

\item[(iii)] {\bf {\mkr Coarse-grained} Kolmogorov operator.} Given a projector $h=\Pi_V$ onto $V$, another useful object related to the conditional expectation \eqref{Eq_G_Markov} is the following coarse-grained Kolmogorov operator: 
\be\label{Eq_reduced_Kolmo}
\overline{\mathcal{K}} f(v)=\int_{\cX_v} \mathcal{K}[f\circ h] (x) \d \mu_v(x),\;\;v\in V, \; \;  f\in C_0^{2}(\V),
\ee
where $\mathcal{K}$ is the Kolmogorov operator defined in Appendix \ref{Sec_ergodic+mixing} for the full SDE.

This operator has an interesting interpretation. Although it is not the ``generator''  of the 
non-Markovian stochastic process $v_t$ associated with the family $(T_t)_{t\geq 0}$, it provides the generator of the Markov process describing the average behavior of $v_t$, when the averaging is taken over the  unobserved variables, i.e.~over the variables lying outside of $V$. 

 \item[(iv)] For the sake of clarity, the formulas \eqref{Eq_G}-\eqref{coeff_G} and \eqref{Eq_G_Markov}-\eqref{Eq_diffusion_averaged} in Theorem \ref{Thm_assymp_gene} have been articulated for the case where $h$ is a projection. Of course such formulas can be generalized to more general mappings $h$, and for $V$ that is not necessarily a subspace of $\mathbb{R}^d$. In such cases, the corresponding coarse-grained Kolmogorov operator, $\overline{\mathcal{K}}$, becomes typically a non-local operator.  Such considerations about the choice of reduced state space has its importance in practice. Indeed non-local features  of $\overline{\mathcal{K}}$ can intuitively help keep a ``trace'' in the reduced state space of certain interactions between the observed and unobserved variables that would be otherwise averaged out by the conditional expectation by using standard projections. 
At the same time, a ``bad'' choice of the observable $h$ can lead to a very poor domain $D(G_t)$ in Lemma \ref{Lem_pseudo}. We refer to \cite{baladi2017linear} for an illuminating discussion regarding the impact of the choice of observable (and their regularity) in the context of linear response theory for (deterministic) Axiom A diffeomorphisms. 
\ei 
\er

\subsection{{\mkr Empirical reduced RP resonances and reconstruction of correlations}}\label{Sec_reducedRP}

As pointed out above, RP resonances although useful to describe power spectra and correlation functions, are difficult 
to estimate in practice when the dimension of the state space gets large. 
In practice, only partial observations of the solutions to Eq.~\eqref{Eq_SDE} are available, e.g.~few solution's components. Theorem \ref{Theorem_A_stoch} shows that from partial observations of a system  that lie within a reduced state space $V$ and whose transitions are collected at a time-lag $\tau$ ($t=\tau$ in \eqref{Eq_main}), a (reduced) Markov operator $T_\tau$ with state space $V$ can be inferred from these observations. Theorem \ref{Theorem_A_stoch} shows then that this operator $T_\tau$ characterizes the coarse-graining in $V$ of the 
transition probabilities in the full state space.  At an intuitive level if a dominant recurrent behavior occurs within an irregular background, then $T_\tau$ must still ``feel'' this recurrent behavior within $V$, in case this dominant behavior is reflected in $V$.  As pointed out already in~\cite{Chek_al14_RP} such a recurrent behavior is manifested by eigenvalues of $T_\tau$ distributed evenly along an inner circle typically close to the unit circle, or by forming a parabola in the complex plane depending on the representation adopted; see Part II \cite{PartII} and Part III  \cite{PartIII} .

An issue though is that in general $T_\tau$ does not come from a Markov semigroup (acting on functions of the reduced state space $V$), and thus one cannot {\it per se} rely on the theory of RP resonances presented in  Sec.~\ref{Sec_spec_Pt} to decompose e.g.~correlation functions of $V$; see Remark \ref{Rmk_spectral-issue}.  Theorem \ref{Thm_assymp_gene} on the other hand, shows, roughly speaking, that $T_\tau$ has a pseudo-generator $G_\tau$ that is close (in a weak sense) to the generator $\mathcal{G}$ of the conditional expectation \eqref{Eq_G_Markov}, when either $\tau$ is sufficiently small or large; see Remarks \ref{Rmk_approx_RPs}-(ii). Thus, in such cases, the resonances estimated from data whose transitions are collected 
at a time-lag $\tau$ in the reduced state space have an interesting interpretation. They approximate  the (weak) RP resonances associated with the generator $\mathcal{G}$ of  the reduced system \eqref{Eq_G_Markov} (Corollary \ref{Cor_approx_RP}), and 
are called the {\bf empirical reduced RP resonances}. Two factors play an important role in this approximation: (a) the partition of $V$ used to approximate the Markov operator $T_\tau$ by a transition matrix $\Gamma_\tau$, and (b) the amount of data used for the estimation of the entries of $\Gamma_\tau$ (see \eqref{P_estimator} below); see \cite{crommelin2006fitting,crommelin2009data}.  When the latter is sufficiently large we may infer that the dominant part of the spectrum of  $\Gamma_\tau$ provides a good approximation of that of 
$T_\tau$, for a sufficiently fine partition of $V$.\footnote{We refer to \cite{schutte2013metastability} for useful error bounds regarding the dominant eigenvalues  for certain types of coarse-graining maps.}  Another important factor in this approximation is, as mentioned above, related to the time-lag $\tau$ as it impacts whether the pseudo-generator $G_\tau$ is a good (weak) approximation of the generator $\mathcal{G}$ of the conditional expectation. If all these factors are favorable, the {\mkr empirical reduced RP resonances}, as good approximation of the generator $\mathcal{G}$'s RP resonances,    describe the solution's variability captured by the conditional expectation \eqref{Eq_G_Markov}; Corollary \ref{Cor_decomp} applied to the Markov semigroup generated by Eq.~\eqref{Eq_G_Markov}. These estimated reduced RP resonances inform us thus about the spectral elements of some coarse-grained version of the original SDE generator, but does not provide in general a sharp approximation of the RP resonances of the full generator $K$.
We refer though to \cite{crommelin2011diffusion} for estimates in the context of homogenized diffusion processes. Translated in the language of this article, it is shown there that the reduced RP resonances  (associated with the homogenized diffusion operator) are asymptotically close to the RP resonances of the underlying multiscale diffusion operator, when the time-scale separation is sufficiently large.

In practice, the dimension of $V$ is kept low so that $T_\tau$ can be efficiently estimated via a maximum likelihood estimator (MLE). 
Note that the reduced state space $V$ and $h$ should be also chosen such that
the observed dynamics in $V$ via $h$ carry relevant information on e.g.~the variability of interest contained in the original system's solutions. 

We detail below our estimation procedure of reduced RP resonances. 
First a bounded domain $\mathcal{D}$ of $V$ should be chosen large enough so that ``most realizations'' of the stochastic process $X_t$ solving Eq.~\eqref{Eq_SDE} fall inside $\mathcal{D}$ after application of the observable $h: \mathbb{R}^d \rightarrow V$,  i.e. $\mathcal{D}$ must be chosen so that $h(X_t)$ belongs to $\mathcal{D}$ for many realizations of the noise in Eq.~\eqref{Eq_SDE}.  This domain is then discretized as the union of $M$ disjoint boxes $B_j$, forming thus a partition.

We assume that our observations in $V$ are made out of the stochastic process $X_t$ (solving Eq.~\eqref{Eq_SDE} in $\mathbb{R}^d$) at discrete time instants
$t=t_n$,  given as multiple of a sampling time $\delta t$, i.e.~$t_n=n\delta t$ with $1\leq n\leq N$, with $N$ assumed to be large. We also assume the time-lag $\tau$ to be equal to a multiple of $\delta t$, i.e.~$\tau=\ell \delta t$.
These observations made in the observation space $V$ are denoted by $Y_n=h(X_{t_n})$.  By adapting the material contained in e.g.~the Supporting Information of~\cite{Chek_al14_RP} (see also \cite{tantet2015early}), the Markov operator $T_\tau$ is approximated by the $M\times M$ transition matrix $\Gamma_\tau$ whose entries 
are given 
\begin{equation}\label{P_estimator}
(\Gamma_{\tau})_{ij}=\frac{\#\bigg\{\Big( Y_{n}\in  B_j \Big)\wedge \Big(Y_{n+\ell } \in  B_i \Big)\bigg\}}{\#\Big\{Y_{n} \in  B_j \Big\}},
\end{equation}
where the $B_j$ form a partition (composed of $M$ disjoint boxes) of the aforementioned domain $\mathcal{D}$ in $V$; see also \cite{crommelin2006fitting,schutte2013metastability} and references therein. In~\eqref{P_estimator}, the notation $\#\{(Y_{n}\in  B_k)\}$ gives the number of observations $Y_{n}$ visiting the box $B_k$, and the logical symbol ``$\wedge$'' means ``and." The leading eigenvalues of the transition matrix $\Gamma_\tau$ can then be computed
with an iterative algorithm such as ARPACK \cite{Lehoucq1997}. We mention that the Markov operator $T_\tau$ can be also approximated --- following the route paved by Ulam \cite{ulam1964collection} ---  from many short simulations instead of using long time series like here; see \cite{froyland2001extracting,dellnitz2001algorithms}. Similarly, the Ulam's method can be adapted for the generator itself to estimate,  when the dimension of the problem permits, the generator's spectral elements without trajectory integration, see \cite{froyland2013estimating}. 

The empirical reduced RP resonances are then obtained as the eigenvalues $\lambda_k(\tau)$ obtained 
from the eigenvalues $\zeta_k(\tau)$ of the Markov matrix $T_\tau$, according to
\begin{equation}\label{Eq_lambda}
\lambda_k(\tau)	= \frac{\log \big(| \zeta_k(\tau)|\big)}{\tau} 
	+  i \frac{\arg \big(\zeta_k(\tau)\big)}{\tau}, \qquad 1\leq k\leq M,
\end{equation}
where $\arg(z)$ (resp.~$\log(z)$) denotes the principal value of the argument (that we adopt to lie in $[-\pi, \pi)$ in this article) (resp.~logarithm) of the complex number $z$.
At a basic level, the motivation behind~\eqref{Eq_lambda} is that the eigenvalues of $\Gamma_\tau$ as the eigenvalues of a Markov matrix, lie within the unit circle (representation that was adopted in~\cite{Chek_al14_RP}) whereas we want here to relate these eigenvalues with the RP resonances associated with the generator $K$ of the original Eq.~\eqref{Eq_SDE}.  This way, the $\lambda_k(\tau)$ given by~\eqref{Eq_lambda} lie naturally within the left-half complex plane.

For $\tau$ sufficiently small or large, i.e.~when the generator of the conditional expectation is reasonably well approximated by the pseudo-generator of $T_\tau$ (in a weak sense, see Theorem \ref{Thm_assymp_gene} and Remark \ref{Rmk_approx_RPs}-(ii)), one can thus proceed as follows to 
measure the amount of solution's variability captured by the conditional expectation:
\bi
\item[($\mathbf{C_1}$)\hspace{-1.5ex}] \hspace{1ex} Given an observable and a reduced state space $V$, one forms the right-hand side (RHS) of \eqref{Eq_decomp_corr2}  in which the $\lambda_j$ are replaced by  the empirical reduced RPs, i.e.~$\lambda_j(\tau)$ given by\eqref{Eq_lambda}, and the eigenfunctions 
by those of the Markov operator $T_\tau$. This way, one forms what we call the ``reduced" correlation function.
\item[($\mathbf{C_2}$)\hspace{-1.5ex}] \hspace{1ex}  One compares the resulting reduced correlation function obtained in step ($\mathbf{C_1}$), to the correlation function as estimated from standard techniques.  If 
the correlations are close, one can conclude that the conditional expectation provides actually a good reduced system and that the non-Markovian effects are negligible to obtain an efficient closure in $V$, for this observable. 
\ei

Note that by replacing the RP resonances by the (empirical) reduced RP resonances in  \eqref{eq:spectralPower}, power spectral densities can  also be approximated in the same fashion.  One should however emphasize that such an analysis depends on the good choice of the time-lag $\tau$ at which the data are collected \cite{tantet2018crisis}, and one should keep in mind this important aspect before drawing any conclusion. 

When $\tau$ is neither small or large, the precise relationships between the $\lambda_k(\tau)$ and the actual RP resonances are non-trivial to characterize in general. Nevertheless, in certain cases, as shown in Part III \cite{PartIII},  the reduced RP resonances are very useful to diagnose and characterize important dynamical features such as nonlinear oscillations embedded within a stochastic background.   
The next section provides an example for which the non-Markovian effects are negligible to capture the slow variable's variability (in terms of correlation functions), whereas the time-scale separation between the ``slow'' and ``fast'' variables is not large.   

\section{Applications to a stochastic slow-fast system}\label{Sec_Langevin}
\subsection{The model}
We consider the following stochastic system
\begin{subnumcases}{\label{Eq_3DHopf}}
\d x =(\lambda x- f y -\gamma xz) \d t + \sigma \d W^1_t\\
\d y =(f x+\lambda y -\gamma yz) \d t+ \sigma \d W^2_t\\
\d z =-\frac{1}{\epsilon}(z-x^2-y^2)\d t+ \frac{\sigma}{\sqrt{\epsilon}} \d W^3_t \label{Eq_3DHopf_z}.
\end{subnumcases}

The stochastic processes $W_t^1, W_t^2, W_t^3$ are independent Brownian motions.
The parameters $\lambda, f, \gamma$ and $\epsilon$ are specified below. In the case $\sigma = 0$, system \eqref{Eq_3DHopf} arises in fluid dynamics and is investigated in \cite{noack_al_2003} as low-dimensional reduced model for a flow past a circular cylinder. 
Hereafter, we deal exclusively with the stochastic case, i.e.~when $\sigma>0$. The reduction problem of this system is analyzed rigorously in \cite{chekroun2019grisanov} using a different approach inspired from \cite{CLW15_vol2}.
In particular it is proved in \cite{chekroun2019grisanov} that this system generates a Markov semigroup which is strong Feller and irreducible that possesses a unique ergodic  invariant measure $\mu$ and thus the theory of RP resonances presented in this article applies.   

System \eqref{Eq_3DHopf} is a slow-fast system driven by additive noise. The theory of slow-fast systems in the deterministic case, i.e.~when $\sigma = 0$ in system \eqref{Eq_3DHopf}, is  well established when the time-scale separation is strong; see for instance  \cite{Jones1995,nipp2013invariant} or the recent monograph \cite{Kuehn:2015tol} and references therein. 
 Indeed, a typical behavior of such systems is characterized by a separation of time scales between the so-called ``slow" and ``fast" variables,as controlled by the parameter $\epsilon$. For our present system, as $\epsilon$ gets small, the $z$-variable exhibits fast fluctuations on timescales over which the $x$- and $y$-variables vary more and more slowly.

Central to the study of such systems in the deterministic context, is the existence of a \textit{slow manifold} which expresses often a (possibly approximate) slaving relationships between the slow- and fast-variables; the latter lying typically in an ``$\epsilon$-neighborhood'' (for $\epsilon$ sufficiently small) to the graph of a function $h$ of the slow variables. In the case of system \eqref{Eq_3DHopf} the slow manifold is explicitly given, for $\sigma=0$, by
\bea \label{Eq_slow_mnf_int}
\Phi: \, & \mathbb{R}\times \mathbb{R} \longrightarrow \mathbb{R},\\
&(x,y) \mapsto x^2 + y^2.
\eea
More precisely,  for $\sigma=0$, any solution $(x(t),y(t),z(t))$ to \eqref{Eq_3DHopf} is attracted exponentially fast to the manifold given by $\mathcal{M}_\epsilon=\mbox{graph} (\Phi) + \mathcal{O}(\epsilon)$, and in particular  the near to slaving relationship $z(t)=h(x^2(t)+y^2(t)) + \mathcal{O}(\epsilon)$ holds for $t$ sufficiently large. Foundations of such a dynamical behavior are found in \cite{FENICHEL197953,10.2307/24890380} for much more general systems than \eqref{Eq_3DHopf}, and we refer to \cite[Chap.~10]{nipp2013invariant} for an introduction to the theory of invariant manifolds (IMs) for singularly perturbed ordinary differential equations (ODEs); see also \cite{Kuehn:2015tol}.

In what follows we determine the empirical reduced RP resonances as described in Sec.~\ref{Sec_reducedRP} not only for the system \eqref{Eq_3DHopf}, but also 
for the reduced system based on the slow manifold $\Phi$, namely by integrating 
\bea\label{Eq_2DHopf}
\d u&=\big(\lambda u- f  v -\gamma  u ( u^2 +  v^2) \big)\d t + \sigma \d W^1_t\\
\d v &=\big(f  v +\lambda  v -\gamma  v (u^2 +  v^2) \big)\d t+ \sigma \d W^2_t.
\eea
The RP resonances of this system, a Hopf normal form subject to an additive white noise, are analyzed in details in Part II \cite{PartII}. 
This system is a ``{\it paragon}'' of nonlinear oscillations in presence of noise.   

The reduced state space $V$ is taken to be the $(x,y)$-plane, being thus only a reduced state space for the original system.
As shown below, the estimation of (reduced) RP resonances allow also for comparing the original system with its slow manifold reduction \eqref{Eq_2DHopf}.

\subsection{Numerical results}

Both systems are numerically integrated via an Euler-Maruyama scheme with a time step of size $\delta t =10^{-5}$.
In each experiment, the systems are simulated up to $T=8\times10^4$ after removal of a transient dynamics of length $T_s=10^3$. 
With $t_n=n \delta t$, the variables collected are $Y_n=(x_n,y_n)$ where $x_n=x(t_n)$ and $y_n=y(t_n)$ for system  \eqref{Eq_3DHopf}, on one hand, and $Y_n=(u_n,v_n)$  where $u_n=u(t_n)$ and $v_n=v(t_n)$ for system \eqref{Eq_2DHopf}, on the other. 
The domain $\mathcal{D}$ used for estimating the Markov matrix $\Gamma_\tau$ in \eqref{P_estimator}, is taken to be 
$\mathcal{D}=[-6, 6] \times [-6, 6]$, decomposed into a uniform grid constituted of $300 \times 300$ cells. 

Three parameter regimes are considered hereafter: two regimes with a strong time-scale separation (Cases I and II), and one with no apparent time-scale separation (Case III); see 
Tables \ref{Table_CaseI} and \ref{Table_CaseIII}.
The time-lag $\tau$ to estimate the transitions is chosen depending on these regimes as follows: $\tau=10^{-3}$ for Case I, and $\tau=10^{-2}$ for Cases II and III.  

\begin{table}[h] 
	\caption{Parameter regimes: Case I and Case II}
	\label{Table_CaseI}
	\centering
	\begin{tabular}{llllll}
		\toprule\noalign{\smallskip}
		& $\lambda$ & $f$ &  $\gamma$  & $\epsilon$ & $\sigma$\\ 
		\noalign{\smallskip}\hline\noalign{\smallskip}
		Case I  & $10^{-3}$  & $10^2$ & $5.6 \times 10^{-2}$ & $10^{-2}$ & $0.55$ \\ 
		Case II & $10^{-3}$ &  $10$ & $1$ & $10^{-2}$ & $0.2$\\  
		\noalign{\smallskip} \bottomrule 
	\end{tabular}
\end{table}

\begin{table}[h] 
	\caption{Parameter regime: Case III}
	\label{Table_CaseIII}
	\centering
	\begin{tabular}{llllll}
		\toprule\noalign{\smallskip}
		& $\lambda$ & $f$ &  $\gamma$  & $\epsilon$ & $\sigma$\\ 
		\noalign{\smallskip}\hline\noalign{\smallskip}
		Case IV  & $10^{-3}$  & $10$ & $1$ & $10$ & $0.3$ \\ 
				\noalign{\smallskip} \bottomrule 
	\end{tabular}
\end{table}

\begin{figure}[htbp]
	\begin{subfigure}{0.48\textwidth}
		\includegraphics[width=1\textwidth,height=0.7\textwidth]{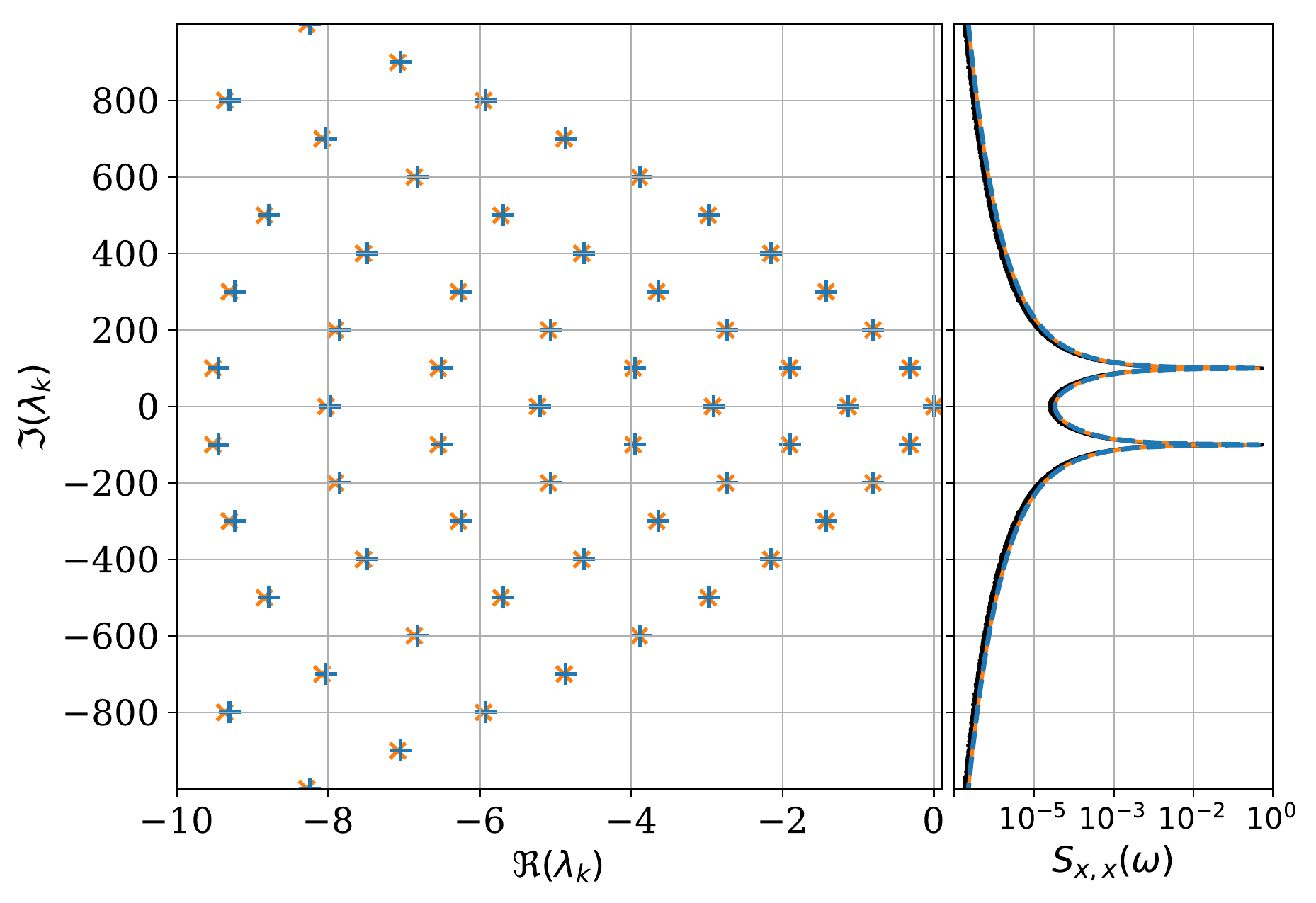}
		\caption{Reduced RP resonances and PSD for Case I}
	\end{subfigure}
	\begin{subfigure}{0.48\textwidth}
		\		\includegraphics[width=1\textwidth,height=0.7\textwidth]{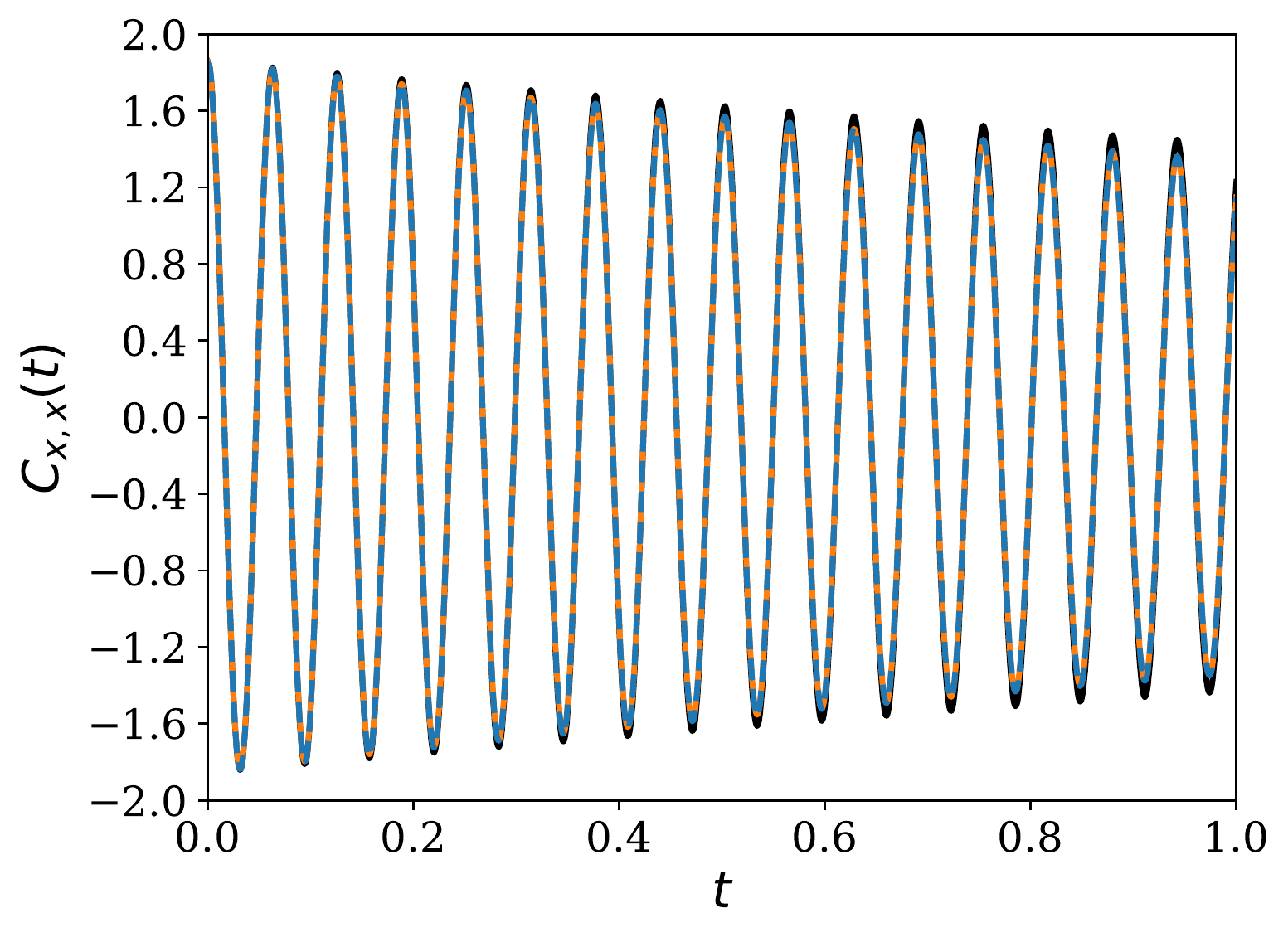}
		\caption{ACF for Case I}
	\end{subfigure}
	\caption{{\footnotesize 
{\bf Leftmost frame in Panel (a):} RP resonances of \eqref{Eq_2DHopf} (orange 'x' sign) and reduced RP resonances of  \eqref{Eq_3DHopf} (blue '+' sign) in the $(x,y)$-plane.
 {\bf Rightmost frame in Panel (a):}  Power spectral densities (PSDs). {\bf Panel (b):} Autocorrelation function (ACF) of $x(t)$ (resp.~$u(t)$).
The spectral reconstructions of the PSD and ACF of $u(t)$ (resp.~$x(t)$) based on  RP resonances (resp.~reduced RP resonances) according to procedure ($\mathbf{C_1}$) described in Sec.~\ref{Sec_reducedRP}, are shown by a dashed orange (resp.~blue) curve, while the black curves show their standard sample estimates.}} 
	\label{Viz_1}
\end{figure}

\begin{figure}[htbp]
	\begin{subfigure}{0.48\textwidth}
		\includegraphics[width=1\textwidth,height=0.7\textwidth]{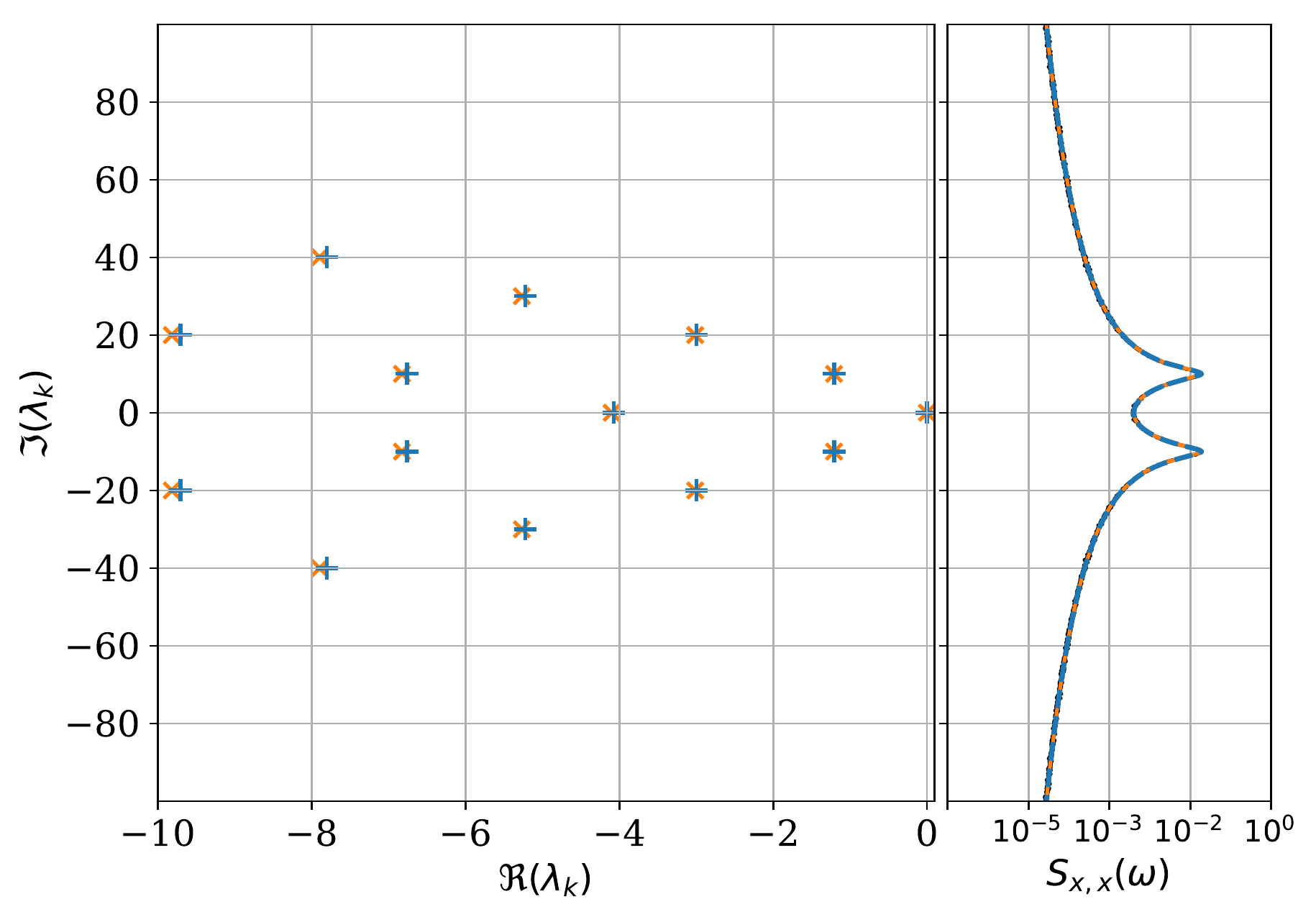}
		\caption{Reduced RP resonances and PSD for Case II}
	\end{subfigure}
	\begin{subfigure}{0.48\textwidth}
		\		\includegraphics[width=1\textwidth,height=0.7\textwidth]{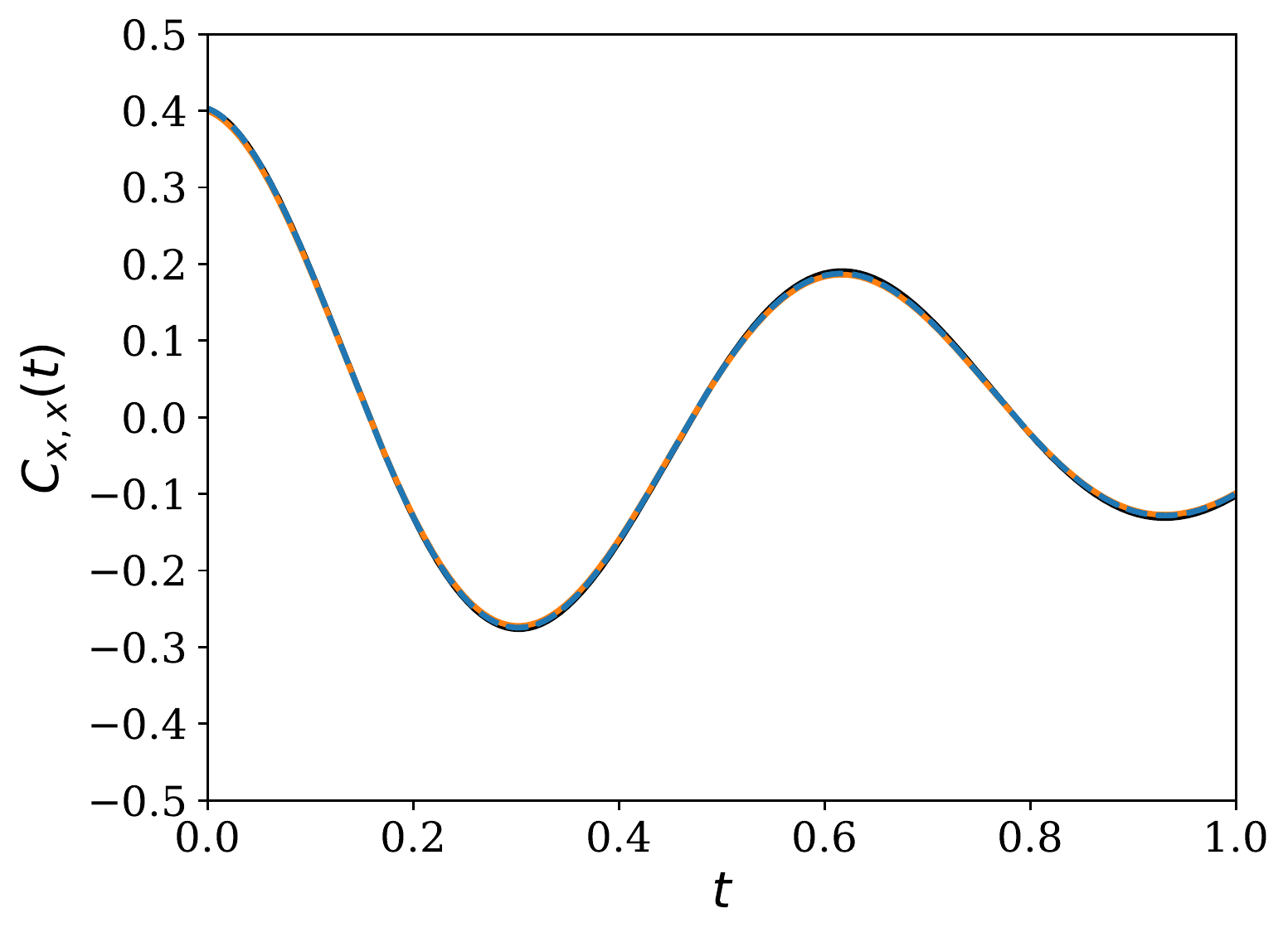}
		\caption{ACF for Case II}
	\end{subfigure}
	\caption{Same as in Fig.~\ref{Viz_1} but for Case II.} 
	\label{Viz_2}
\end{figure}

\begin{figure}[htbp]
	\begin{subfigure}{0.48\textwidth}
		\includegraphics[width=1\textwidth,height=0.7\textwidth]{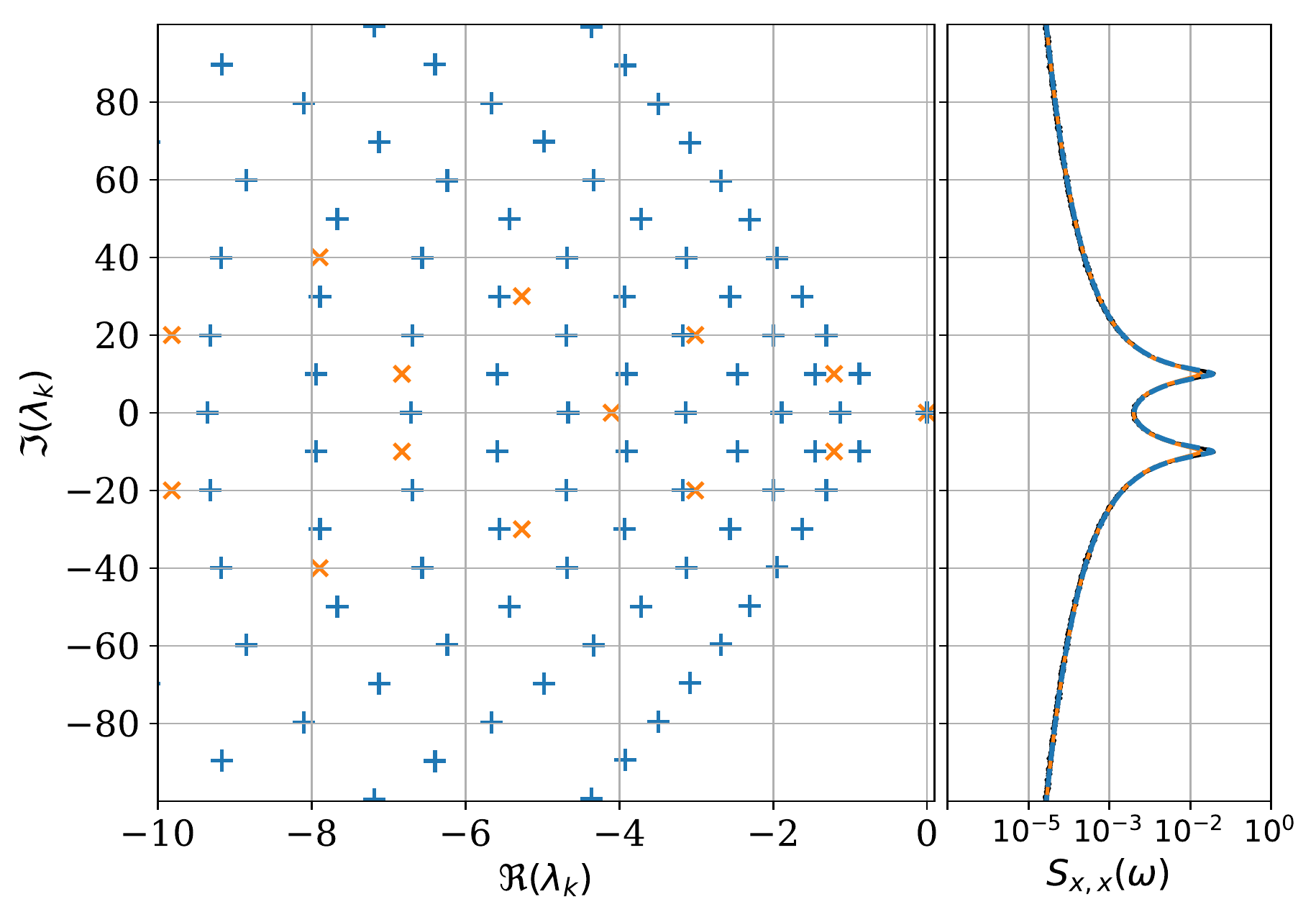}
		\caption{Reduced RP resonances and PSD for Case III}
	\end{subfigure}
	\begin{subfigure}{0.48\textwidth}
		\		\includegraphics[width=1\textwidth,height=0.7\textwidth]{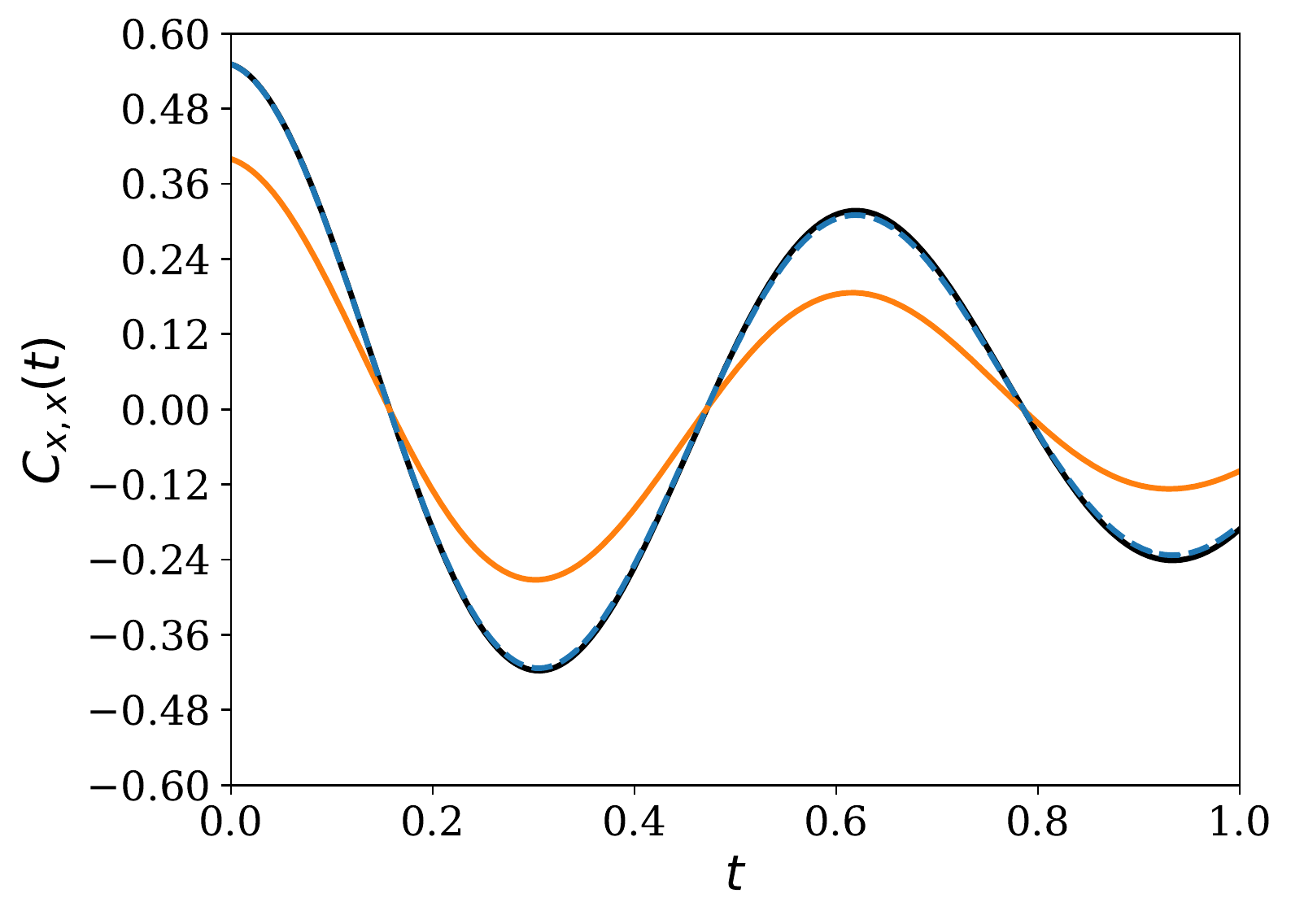}
		\caption{ACF for Case III}
	\end{subfigure}
	\caption{Same as in Fig.~\ref{Viz_1} but for Case III.} 
	\label{Viz_3}
\end{figure}

Within this experimental protocol, we first estimate the reduced RP resonances.  To do so, we first estimate the Markov matrix entries of $\Gamma_\tau$ according to \eqref{P_estimator}, and then determine the reduced RP resonances according to \eqref{Eq_lambda}. For the three parameter regimes considered here, these resonances 
are shown by blue '+' signs in each Panel (a) of Figures \ref{Viz_1}, \ref{Viz_2}, and \ref{Viz_3}.   
The resonances for the slow manifold system \eqref{Eq_2DHopf} are estimated according the same procedure. They are no longer reduced but genuine 
RP resonances since the system is two-dimensional.  Thus, the RP resonances associated with system \eqref{Eq_2DHopf} are shown by orange 'x' signs 
in each Panel (a) of Figures \ref{Viz_1}, \ref{Viz_2}, and \ref{Viz_3}.   
In each of these Panels (a),  the rightmost narrow panel shows the spectral reconstructions of the PSD of $u(t)$ (resp.~$x(t)$) based on  RP resonances (resp.~reduced RP resonances) according to  procedure ($\mathbf{C_1}$) described in Sec.~\ref{Sec_reducedRP}, by a dashed orange (resp.~blue) curve, while the black curve shows its standard PSD sample estimate.  In each Panel (b) of  Figures \ref{Viz_1}, \ref{Viz_2}, and \ref{Viz_3} are shown the autocorrelation function (ACF) of $u(t)$ (resp.~$x(t)$) based on  RP resonances (resp.~reduced RP resonances) based on the procedure ($\mathbf{C_1}$), and here again the black curve its standard ACF sample estimate.  

The conclusions of these numerical experiments are without ambiguity. First, in the case of a strong  time-scale separation between the observed variables ($x$,$y$) and the unobserved variable ($z$), i.e.~in Cases I and II,  the reduced RP resonances allow for an almost exact reconstruction of the ACF and PSD of $x(t)$ (as for $y(t)$, not shown).
This is explained from the theoretical understanding provided by  Secns.~\ref{Sec_cond_exp} and \ref{Sec_reducedRP}. Indeed, as explained therein, since the time-lag $\tau$ is small here, it is expected that the reduced RP resonances provide a good approximation of the resonances associated with the generator of the conditional expectation \eqref{Eq_G_Markov}.  On the other hand,  the conditional expectation is known to 
provide a (very) good approximation of the dynamics of the slow variables when $\epsilon$ is small, and that it coincides with the slow manifold reduced system \eqref{Eq_2DHopf} as $\epsilon \rightarrow 0$; see \cite{chekroun2019grisanov}. Thus without any surprise the RP resonances of  system \eqref{Eq_2DHopf} coincide with the reduced RP resonances of  \eqref{Eq_3DHopf}  (associated with the oberved variables $(x,y)$) for a 
strong  time-scale separation. The fact that these resonance are organized  along parabolas in the left half plane, is thoroughly analyzed in Part II \cite{PartII}.

The conclusion for the case of no time-scale separation (Case III) is more subtle, but reveals an interesting usefulness of the RP resonance analysis.  
For this case, the RP resonances of the slow manifold reduced system \eqref{Eq_2DHopf} differ from the reduced RP resonances (see Panel (a) in Fig.~\ref{Viz_3}) which come 
here with a noticeable  mismatch at the level of ACFs although not revealed at the level of PSDs.  On the contrary, the reduced RP resonances allow for an almost perfect reconstruction of the ACF.  As a consequence, one can conclude that the slow manifold is no longer a valid parameterization of the $z$-variable in terms of the $x$- and $y$-variables (also observed in \cite{chekroun2019grisanov} for this regime) but that still, the conditional expectation \eqref{Eq_G_Markov} provides a valid approximation for the reduced dynamics.    

Such a diagnosis regarding the conditional expectation drawn from RP resonance analysis offers promising perspectives in terms of reduced-order modeling. As illustrated here, the latter resonances could indeed serve to determine whether efforts on non-Markovian effects must be invested or not. In that respect, we mention the recent variational approach relying on optimal manifolds minimizing a parameterization defect \cite{CLM19_closure} which offers new perspectives to approximate analytically from the original equations the conditional expectation; see also \cite{CLM16_Lorenz9D}.  Applied to Case III, this approach based on optimal parameterizing manifolds (OPMs) provides an analytical substitute to the slow manifold which allows in turn for the derivation of an efficient 2D reduced system of SDEs (without non-Markovian terms), confirming the conclusions drawn from the RP resonance analysis conducted here.

\begin{acknowledgements}
\vspace{-1ex}
The authors would like to thank the reviewers for their very useful and constructive comments.
This work has been partially supported by the European Research Council under the European Union's Horizon 2020 research and innovation program (grant agreement No.~810370 (MDC)),  by the Office of Naval Research (ONR) Multidisciplinary University Research Initiative (MURI) grant N00014-16-1-2073  (MDC), by the National Science Foundation grants OCE-1658357, DMS-1616981(MDC), AGS-1540518 and AGS-1936810 (JDN), by the LINC project (No.~289447) funded by EC's Marie-Curie ITN (FP7-PEOPLE-2011-ITN) program (AT and HD) and by the Utrecht University Center for Water, Climate and Ecosystems (AT).
\end{acknowledgements}

\newpage
\appendix

\section{Elements of stochastic analysis}\label{Sec_Appendix}
In this appendix we present a short survey of elements of stochastic analysis used in the Main Text. 
The main objective is to introduce the key concepts and tools of stochastic analysis for stochastic differential equations (SDEs), to a wider audience in the geosciences and macroscopic physics. 

\subsection{Markov semigroups}\label{Sec_Marov_smg}

Two approaches dominate the analysis of stochastic dynamics. We are here concerned with the approach rooted in Stochastic Analysis which, contrary to the random dynamical system (RDS) approach \cite{Arnold98,CSG11,CLW15_vol1}, does not substitute a deterministic (nonlinear) flow $S(t)$ by a stochastic flow $S(t,\omega)$ acting\footnote{$\omega$ labelling the noise realization.} on the state space $\cX$ but rather by a family of linear operators $P_t$, acting on a space of observables of the state space, i.e.~on functions of $\cX$.  A typical choice of observables is given by $\mathcal{C}_b(\cX)$, the space of bounded and continuous functions on $\cX$. In what follows $\cX$ is a finite-dimensional Polish space.

More precisely, this family $P_t$ reflects the (averaged) action of the stochastic flow at the level of functions and is  
given as the mapping which to each function $\phi$ in $\mathcal{C}_b(\cX)$ associates the function:
\be\label{Pt_def}
P_t \phi (x) =\mathbb{E}( \phi(S(t,\cdot)x))=\int_{\Omega} \phi (S(t,\omega) x) \d \mathbb{P} (\omega), \quad t\geq 0,\; x\in \cX.
\ee
In \eqref{Pt_def}, the function $\phi$ is the aforementioned {\it observable}. Its physical meaning could be, for instance, the potential vorticity or the temperature of a 
fluid at a given location or averaged over a volume. The RHS of \eqref{Pt_def} involves averaging over the realizations $\omega$, i.e.~expectation.   For deterministic 
flow it reduces to $P_t\phi(x)=\phi(S(t)x)$ and is known as the Koopman operator. Note that $P_t$ such as defined in \eqref{Pt_def} is not limited to stochastic flow, more generally
 $P_t \phi (x) =\mathbb{E}( \phi(X_t^x))$ where  $X_t^x$ denotes a stochastic process  that solves Eq.~\eqref{Eq_SDE} (as associated with $P_t$) and emanates from $x$ in $\cX$.

Under general assumptions on $F$ and $D$, the stochastic process $X_t$ solving Eq.~\eqref{Eq_SDE} is Markovian (i.e.~the future is determined only by the present value of the process) which translates at the level of $P_t$ into the following semigroup property
\be\label{Eq_smg_pty}
P_0=\mbox{Id}, \quad P_t  P_s=P_{t+s}, \; \; t,s\geq 0.
\ee  
A breakdown of \eqref{Eq_smg_pty} indicates thus that the underlying stochastic process is non-Markovian. 

It is noteworthy to mention that even when $P_t$ satisfies \eqref{Eq_smg_pty}, it does not ensure that $P_t$ is a {\it strongly continuous} semigroup \cite{Pazy1983} on $\mathcal{C}_b(\cX)$. Nevertheless, $(P_t)_{t\geq 0}$ is extendable  to a strongly continuous semigroup  in $L^2_{\mu}$ as soon as $\mu$ is an invariant measure of the Markov semigroup; see Theorem \ref{Thm_cond_mixing} below. The spectral theory of such semigroups \cite{Engel_Nagel}  {\mk is at the core of the} description of mixing properties in $L^2_{\mu}$,  such as presented in Sec.~\ref{Sec_spec_Pt} in the Main Text. 

\subsection{Ergodic invariant measures and the strong Feller-Irreducibility approach}\label{Sec_Feller}

The Fokker-Planck equation \eqref{Eq_FKP} may support several weak stationary solutions. An important question, is thus the identification of stationary measures that describe the asymptotic statistical behavior of the solutions of Eq.~\eqref{Eq_SDE}, in a typical fashion. The notion of ergodic invariant measures plays a central role in that respect, and relies on the following important characterization of ergodic measures for (stochastically continuous) Markov semigroups \cite[Theorem~3.2.4]{DaPrato1996}.
\bd\label{Def_ergodic}
An invariant measure is ergodic if one of the following three equivalent statements holds:
\bi
\item[i)] For any $f \in L^2_\mu(\cX)$, if $P_t f =f$, almost surely w.r.t $\mu$~($\mu$-a.s.)~for all $t\geq 0$, then $f$ is constant $\mu$-a.s. 
\item[ii)] For any Borel set $\Gamma$ of $\cX$, if $P_t \mathbbm{1}_{\Gamma} =\mathbbm{1}_{\Gamma}$ $\mu$-a.s.~for all $t\geq 0$, then $\mu(\Gamma)=0$ or $1$.
\item[iii)] For any $f \in L^2_\mu (\cX)$, $\frac{1}{T} \int_0^T P_s f \d s \underset{T\rightarrow \infty}\longrightarrow \int f \d \mu$ in $L^2_\mu(\cX).$
\ei
\ed

In practice, an efficient approach to show the existence of an ergodic measure consists of showing the existence of a {\it unique} invariant measure, since in this case such an invariant measure is necessarily ergodic \cite[Theorem~3.2.6]{DaPrato1996}. Various powerful approaches exist to deal with the existence of a unique invariant measure.  The next section discusses the classical approach based on the theory of {\it strong Feller Markov semigroups} and {\it irreducibility}.

The main interest of the strong Feller-Irreducibility approach lies in its usefulness for checking the conditions of the  
Doob-Khasminskii Theorem \cite{doob1948asymptotic,khas1960ergodic,DaPrato1996}, the latter ensuring  the existence of {\it at most} one ergodic invariant measure. This strategy requires the proof of certain smoothing properties of the
associated Markov semigroup, and to show that any point can be (in probability) reached at any time instant
by the process regardless of initial data. This property is known as {\it irreducibility}.  It means that $P_t \mathbbm{1}_{U} (x)> 0$ for all $x$ in $\cX$, every $t>0$, and all non-empty open sets $U$ of $\cX$, which is equivalent to say that 
\be
\mathbb{P}(\|S(t,\cdot)x-z\| < \epsilon) >0, 
\ee
for any $z$ in $\cX$, $\epsilon >0$ and $t>0$; see \cite[p.~67]{cerrai2001second}. In other words the irreducibility condition expresses the idea that any neighborhood of any point $z$  in $\cX$, is reachable at each time, with a positive probability. 

Remarkably,  the irreducibility is usually inferred  from the controllability of the associated control system $\dot{x}=F(x) + D(X) u(t)$; see \cite{chekroun2019grisanov} for a simple illustration. This approach is well-known and based on the support theorem of Stroock and Varadhan \cite{stroock1972} (see also \cite[Theorem 8.1]{ikeda2014stochastic}) that shows that  several properties of the SDEs can be studied and expressed in terms of the control theory of ordinary differential equations (ODEs); see \cite[Secns.~7.3 and 7.4]{DaPrato1996} for the case of additive (non-degenerate) noise and  \cite{arnold1987unique,kliemann1987recurrence} for the more general case of nonlinear degenerate noise, i.e.~in the case where the noise acts only on part of the system's equations, corresponding to ker$(Q)\neq \{0\}$. 

The strong Feller property means that the Markov semigroup maps bounded measurable
functions into bounded continuous functions. This property, related to a regularizing effect
of the Markov semigroup $(P_t)_{t\geq 0}$, is a consequence of the {\it hypoellipticity}
of the {\it Kolmogrorov operator} $\mathcal{K}$ defined on smooth functions $\psi$ (of class $C^2$) as follows when $\cX=\mathbb{R}^d$:
\be\label{Def_K0-additive}
\mathcal{K}\psi(x)= \frac{1}{2}\mbox{Tr}(\mathbf{\Sigma} \nabla^2 \psi(x)) +\langle F(x), \nabla \psi \rangle,
\ee
where 
\be
\mbox{Tr}(\mathbf{\Sigma} \nabla^2 \psi(\cdot))=\sum_{i,j=1}^d \big[D(x)D(x)^T\big]_{ij}\partial_{ij}^2 \psi.
\ee 
Here $\mbox{Tr}$ denotes the trace of a matrix. Note that hypoelliptic operators include those that are uniformly elliptic for which the Weyl's smoothing lemma applies; e.g.~\cite[Theorem~4.7]{dacorogna2004introduction}. Hypoellipticity allows nevertheless for dealing with the  case  of degenerate noise, which is important in applications. 

A very efficient criteria for hypoellipticity is given by  H\"ormander's theorem \cite{hormander1967hypoelliptic,norris1986simplified}; see also \cite[Appendix C1]{CSG11} for a discussion on the related H\"ormander's bracket condition and its implications to the existence of other types of meaningful measures for SDEs, namely the  Sina\"i-Ruelle-Bowen (SRB) random measures.  We refer also to Part II \cite{PartII}, for an instructive verification of the H\"ormander's condition in the case of the Hopf normal form subject to additive noise.

From a geophysical perspective, it is noteworthy to mention that the strong Feller-Irreducibility approach allows for dealing with a broad class of truncations of fluid dynamics models that would be perturbed by noise, possibly degenerate. For instance, in the case of truncations of 2D or 3D Navier-Stokes equations,
 the strong Feller-Irreducibility approach has been shown to be applicable even for an additive noise that forces only very few modes \cite{agrachev2005navier,romito2004ergodicity}. The delicate point of the analysis is the verification of the controllability (and thus irreducibility) of the associated control  system, by techniques typically adapted from \cite{jurdjevic1985polynomial} or rooted in chronological calculus  as in \cite{agrachev2005navier}.  Whatever the approach, the analysis requires the appropriate translation into geometrical terms of the cascade of energy in which the nonlinear terms transmit the forcing from the few modes to all the others \cite{romito2005geometric}.
We mentioned however \cite{majda2016ergodicity} for an example of a stochastic dynamical system which has the square of the Euclidean norm as the Lyapunov function, is hypoelliptic with nonzero noise forcing, and that yet fails to be reachable or ergodic.

\subsection{Markov semigroups and mixing}\label{Sec_ergodic+mixing}

We recall here standard results about Markov semigroups. It states that any Markov semigroup that is strong Feller and irreducible and for which an invariant measure exists (which is thus unique)  is not only ergodic but also strongly mixing for the total variation norm of measures.  
Given two  probability measures $\mu_1$ and $\mu_2$ on $\cX$, we recall that the latter is defined as \cite[Eq.~(3.1)]{hairer2009introduction}   
\be\label{mix_1b}
{\mbox TV}(\mu_1,\mu_2)=\sup_{\substack{ g \in \mathcal{B}_b(\cX)\\ \| g\|_{\infty}\leq 1}} \bigg|  \int g \, \d \mu_1-\int_{\cX}  g \, \d \mu_2\bigg|,
\ee
where $\mathcal{B}_b(\cX)$ denotes the set of Borel measurable and bounded functions on $\cX.$
\begin{thm}\label{Thm_cond_mixing}
Let $\mu$ be an invariant measure of a Markov semigroup $(P_t)_{t\geq 0}$. For any $p\geq 1$ and $t\geq 0$, $P_t$ is extendable to a linear bounded operator on $L^p_\mu(\cX)$ still denoted by $P_t$. Moreover

\bi
\item[(i)] $\|P_t\|_{\mathcal{L}(L^p_\mu(\cX))} \leq 1$
\item[(ii)] $P_t$ is strongly continuous semigroup in $L^p_\mu(\cX)$.
\ei

If furthermore $(P_t)_{t\geq 0}$ is strong Feller and irreducible, then $\mu$ is ergodic (and unique)  and for any $x$ in $\cX$ and $g$  in $L^1_{\mu}$
 \be\label{strong_ergo}
\underset{T\rightarrow \infty}\lim \frac{1}{T} \int_0^T g(X_\tau^x) \d \tau= \int_{\cX} g(x)\d \mu, \; \; \mathbb{P}\mbox{-a.s.},
\ee
where $X_t^x$ denotes the stochastic process solving the SDE associated with $P_t$. 

{\mk In this case, the invariant measure} $\mu$ is also {\bf strongly mixing} in the sense that for any measure $\nu$ on $\cX$, we have:
\be\label{mix_1}
{\mbox TV}(\cL_t \nu ,\mu)\underset{t\rightarrow \infty}\longrightarrow 0.
\ee

\end{thm}

For the definition of a strongly continuous semigroup also known as $C_0$-semigroup we refer to \cite[p.~36]{Engel_Nagel}. 
For an introduction to semigroup theory we refer to \cite{van2012asymptotic,engel2006short}.

\bp[Proof.]
We prove first $(i)$. {\mk The proof is standard and can be found e.g.~in \cite[Prop.~1.14]{guionnet2003lectures} but is reproduced here for the reader's convenience}. Let $g$  be in $\mathcal{C}_b(\cX)$. By the H\"older inequality, we have
\be\label{Ineq_base}
|P_t g (x) |^p \leq  P_t (|g|^p) (x).
\ee

If we now integrate both sides of this inequality with respect to $\mu$, we obtain 
\be\label{Eq_Holder}
\int_{\cX} |P_t g (x) |^p  \mu (\d x) \leq \int_{\cX} P_t (|g |^p) (x) \mu (\d x)= \int_{\cX} | g|^p (x) \mu (\d x),
\ee
the latter equality resulting from the invariance of $\mu$.
Since $\mathcal{C}_b(\cX)$ is dense in $L^p_\mu(\cX)$,  the inequality \eqref{Eq_Holder} can be extended to any function in $L^p_\mu(\cX)$, and thus $(P_t)_{t\geq 0}$ can be uniquely extended to a contraction semigroup in  $L^p_\mu(\cX)$, and property $(i)$ is proved. 

Let us show now that $(P_t)_{t\geq 0}$ is strongly continuous in $L^p_\mu(\cX)$.  
Since $(P_t)_{t\geq 0}$ is a Markov semigroup, for any $g$ in $\mathcal{C}_b(\cX)$ and $x$ in $\cX$, we have that the mapping 
$t \mapsto P_t g (x)$ is continuous. Therefore  by the dominated convergence theorem
\be
\lim_{t\rightarrow 0} \; P_t g  =g  \mbox{ in  } L^p_\mu(\cX).
\ee
The density of $\mathcal{C}_b(\cX)$  in $L^p_\mu(\cX)$ allows us to conclude that this convergence holds 
when $g$ is in $L^p_\mu(\cX)$.

The ergodicity of $\mu$ results from the aforementioned Doob's theorem. The time-average property \eqref{strong_ergo} and the mixing property \eqref{mix_1} can be obtained as a consequence of e.g.~\cite[Cor.~2.3]{seidler1997ergodic}; see also \cite[Cor.~1]{stettner1994}. \qed
\ep

\subsection{Generator of a Markov semigroup}\label{Sec_gen_Markov}
Recall that the generator $A$ of any strongly continuous semigroup $(T(t))_{t\geq 0}$ on a Hilbert space $\mathcal{H}$ is defined as
the operator $A:D(A)\subset \mathcal{H} \rightarrow \mathcal{H}$, such that 
\be
A \varphi = \lim_{t\rightarrow 0^+} \frac{1}{t}\big(T(t) \varphi -\varphi\big),
\ee 
defined for every $\varphi$ in the domain
\be\label{Def_domA}
D(A)=\{\varphi \in \mathcal{H} \;|\; \lim_{t\rightarrow 0^+} \frac{1}{t}\big(T(t) \varphi -\varphi\big) \mbox{ exists}\}.
\ee

As any generator of a contraction semigroup, given an invariant measure $\mu$, the generator $K$ of the contraction semigroup $(P_t)_{t\geq 0}$  in  $L^2_{\mu}$ (Theorem~\ref{Thm_cond_mixing}-(i)) is dissipative, which is equivalent to say, since $L^2_{\mu}$ is a Hilbert space, that 
\be\label{Eq_dissip_K}
\Re \langle K f,f \rangle_{L^2_{\mu}} \leq 0, \; \forall \, f \in D(K),
\ee
where $D(K)$ denotes the domain of $K$; see e.g.~\cite[Prop.~II.3.23]{Engel_Nagel}.
The domain $D(K)$ is furthermore dense in $L^2_{\mu}$ and  $K$ is a closed operator; see \cite[Cor.~2.5 p.~5]{Pazy1983}. 
The isolated part of the spectrum of $K$ provides the Ruelle-Pollicott resonances; see Sec.~\ref{Sec_spec_Pt}.

\subsection{Return to equilibrium and spectral gap}\label{Sec_spec_gap}

We present here some useful results concerning (i) the exponential return to equilibrium  for strong Feller and irreducible Markov semigroups, and (ii) spectral gap in the spectrum of the Markov semigroup generator $K$; see Theorems \ref{Thm_cve1} and \ref{Thm_cve2} below. Theorem \ref{Thm_cve1} deals  with  semigroups that become {\it quasi-compact} after a finite time, and Theorem  \ref{Thm_cve2}  addresses the exponential $L^2$-convergence and lower bound of the spectral gap.  For Theorem \ref{Thm_cve1}, the approach is based on Lyapunov functions such as formulated in \cite{bellet2006ergodic}. We propose a slightly different presentation for which we provide the main elements of the proof. We refer to \cite{douc2009subgeometric} for an efficient (and beautiful) generalization of such Lyapunov-type criteria allowing for {\it sub-exponential convergence} towards the equilibrium.

Recall that the {\it essential spectral radius} $\mathbf{r}_{ess}(T)$ of a linear bounded operator $T$ on a Banach space $\mathcal{E}$ satisfies \cite[p.~249]{Engel_Nagel} the Hadamard formula
\be\label{Spec_rad}
\mathbf{r}_{ess}(T)=\underset{n\rightarrow \infty}\lim \| T^n\|_{ess}^{1/n},
\ee
where 
\be
\|T\|_{ess}=\inf \Big\{ \|T- \mathcal{C}\|_{\mathcal{L}(\mathcal{E})} \, : \, \mathcal{C} \mbox{ is  a linear and compact operator of } \mathcal{E}\Big\}.
\ee
We have then the following convergence result.
\begin{thm}\label{Thm_cve1}
Let $\mathcal{P}=(P_t)_{t\geq 0}$ be a strong Feller and irreducible Markov semigroup in $L^2_\mu(\mathbb{R}^d)$ ($\cX=\mathbb{R}^d$) generated by an SDE given by Eq.~\eqref{Eq_SDE} for which $F$ and $G$ are  locally Lipschitz. Assume that  there exists a Lyapunov function\footnote{Recall that a $C^2$ function $U$  is called  a Lyapunov function $U(x)\geq1$ and $\lim_{|x|\rightarrow \infty} U(x)=\infty$, ensuring thus that the level sets $\{U\leq \alpha\}$ are compact.} $U$ and a compact set $\mathfrak{A}$ for which there exist $a >0$, $0<\kappa<1$ and $b<\infty$, such that
\begin{subequations}\label{Eq_Lyap_cond}
\begin{align}
\mathcal{K} U & \leq a  U,\label{Eq_Lyap_cond-a}\\
P_{t_0} U &\leq \kappa U + b  \mathbbm{1}_{\mathfrak{A}}, \; \mbox{ for some }  t_0>0,\label{Eq_Lyap_cond-b}
\end{align}
\end{subequations}
where $\mathcal{K}$ is the Kolmogorov differential operator generating the Markov process associated with $\mathcal{P}$. 
 Then for all $t>t_0$, $P_t$ becomes quasi-compact, i.e.
\be\label{Eq_quasi}
\mathbf{r}_{ess}(P_t)\leq \kappa,
\ee
where the essential spectral radius is taken for $P_t$ as acting on  $\mathcal{E}=\mathcal{F}_{U}$ given by
\be\label{Eq_FU}
\mathcal{F}_{U}=\{f:\mathbb{R}^d \rightarrow  \mathbb{R} \;|\;  \mbox { f Borel measurable  and }  \| f\|_{U} <\infty\},
\ee
and endowed with the norm 
\be
\|f\|_U=\underset{x\in \mathbb{R}^d}\sup \frac{|f(x)|}{U(x)}.
\ee 

Furthermore $(P_t)_{t\geq 0}$ has a unique invariant measure $\mu$, and the inequality \eqref{Eq_quasi}  ensures that there exist $C>0$ and $\lambda >0$ such that for all $f$ in $\mathcal{F}_{U}$,
\be\label{CVE_pre_totvar}
\left| P_t f(x) -\int f \d \mu \right| \leq C e^{-\lambda t} U(x), \;\; t>t_0, \; \forall \, x\in \mathbb{R}^d. 
\ee

\end{thm}

{\mk The proof of this result is found in Appendix \ref{Sect_proof}. }  
 \br
 The assumption \eqref{Eq_Lyap_cond-b} {\mk is sometimes} verified from moment estimates in practice.  For instance if there exist $k_0>0$ and $k_1>0$ such that
 \be\label{Eq_toto0}
 \mathbb{E} |X_t^x| \leq k_0 e^{-k_1 t}|x| +c, \; t\geq 0,
 \ee
 then for any $t\geq -\frac{1}{k_1} \log(\frac{1}{4 k_0})$, we have $ \mathbb{E} (|X_t^x| +1)\leq \frac{1}{2} (|x|+1)-\frac{1}{4} |x|+c+\frac{1}{2}$, which leads to
  \be
 \mathbb{E} (|X_t^x| +1) \leq \frac{1}{2} \big(|x|+1\big) +\big(c+\frac{1}{2} \big)\mathbbm{1}_{B_r},
 \ee
 for all $r>4(c+\frac{1}{2}),$ and thus \eqref{Eq_Lyap_cond-b} holds with $U(x)=|x|+1.$

More generally, if 
\be\label{Meta_Lyap_cond}
\mathcal{K} U  \leq -\alpha U +\beta,  \; \mbox{ with }\alpha >0,  \mbox{ and } \; 0\leq \beta<\infty,
\ee
then $\frac{\d}{\d t} P_t U (x)={\mk P_t \mathcal{K} U} (x)\leq -\alpha P_t U(x) +\beta$, leading to 
\be\label{Eq_toto1}
\mathbb{E} \big[U(X_t^x)\big]\leq U(x) e^{-\beta t} +\frac{\beta}{\alpha} \big( 1-e^{-\alpha t}\big), \;\; t>0,
\ee
and similarly \eqref{Eq_Lyap_cond-b} holds.  In addition, \eqref{Meta_Lyap_cond} implies \eqref{Eq_Lyap_cond-a}.
{\mk Note that \eqref{Meta_Lyap_cond} and  \eqref{Eq_toto1} are quite standard; see e.g.~\cite[Lemma 2.11]{dragoni2012ergodicity}.} 
 
{\mk Finally, note also that finding a Lyapunov function may be easier than proving inequalities of the form \eqref{Eq_toto0}. 
For instance, if there is a Lyapunov function which grows polynomially like $\|p\|^q$, then one knows that the process has moments of order $q$; see \cite{meyn1993stability2,meyn1993stability3}. } 
 \er

Finally, {\mk lower bounds of the spectral gap in $L^2_{\mu}$ may be derived for a broad class of 
 SDEs}.  Recall that the generator $K$ has a spectral gap in $L^2_{\mu}$ if there exists $\delta>0$ such that
 \be\label{Cond_gap}
 \sigma(K)\cap \{\lambda\,:\, \mathrm{Re}(\lambda)>-\delta\}=\{0\}.
 \ee
 The largest $\delta>0$ with this property is denoted by $\mathrm{gap}(K)$, namely
 \be\label{Formula_gap}
 \mathrm{gap}(K)=\sup \{ \delta >0 \mbox{ s.t. }\eqref{Cond_gap} \mbox{ holds} \}.
 \ee

The following result is a consequence in finite dimension of more general convergence results \cite[Theorems 2.5 and 2.6]{goldys2005}. Since $(P_t)_{t\geq 0}$ is a C$_0$-semigroup in $L^2_{\mu}$,  the theory of asymptotic behavior of a semigroup with a strictly dominant, algebraically simple eigenvalue (e.g.~\cite[Theorem.~3.6.2]{van2012asymptotic}) implies the spectral gap property stated in the following.  
\begin{thm}\label{Thm_cve2}
Assume that $(P_t)_{t\geq 0}$ is strong Feller and irreducible.  Assume furthermore 
that the following ultimate bound holds for the associated stochastic process $X_t^x$,
i.e.~there exist $c,k,\alpha>0$ such that
\be
\mathbb{E}  \, |X_t^x|^2< k |x|^2e^{- \alpha t} +c, \; \;t\geq 0, \; \; x\in \mathbb{R}^d.
\ee
Then there exists a unique invariant measure $\mu$ for which the $U$-uniform ergodicity \eqref{CVE_pre_totvar} holds with $U(x)=1+|x|^2$,  as well as the following exponential $L^2$-convergence 
 \be\label{CVE_inL2}
 \| P_t \varphi - \int \varphi \d \mu \|_{L^2_{\mu}} \leq C e^{-\lambda t} \|\varphi \|_{L^2_{\mu}}, \; t\geq 0, \; \varphi \in L^2_{\mu},
 \ee
with $C$ and $\lambda$ positive constants independent of $\varphi$; the latter rate of convergence being the same as that of \eqref{CVE_pre_totvar}. Furthermore, one has the following lower bound for the $L^2_{\mu}$-spectrum of the generator $K$:
\be
0<\lambda\leq \mathrm{gap}(K).
\ee 
 \end{thm}
 
We will see in Part II \cite{PartII} of this three-part article that Theorem \ref{Thm_cve2} has important practical consequences. In particular it shows for a  broad class of controllable ODEs, perturbed by a white noise process for which the Kolmogorov operator is hypoelliptic, that an $L_{\mu}^2$-spectral gap is naturally induced by the noise whereas in absence of the latter the gap may be zero, leading thus to a form of mixing enhancement by the noise.  We finally mention \cite{hairer2014spectral} for other conditions, ensuring an $L^2_\mu$-gap based on spectral gaps in Wasserstein distances, verifiable in practice by following the approach of \cite{hairer2011asymptotic}. 

\subsection{Proof of Theorem \ref{Thm_cve1}}\label{Sect_proof}

\bp[Proof.]

It is standard from the theory of Lyapunov functions that the existence of a unique invariant measure $\mu$ is ensured by the condition \eqref{Eq_Lyap_cond-a} together with the irreducibility and strong Feller properties. The rest of the proof is thus concerned with \eqref{Eq_quasi} and the exponential convergence \eqref{CVE_pre_totvar}. 
 
{\bf Step 1.}   First, note that the It\^o formula gives
\be
\d U =  \mathcal{K} U \d t + \mbox{ ``Martingale''},
\ee
which leads  (since $\mathcal{K} U \leq a  U$) to
\be
\mathbb{E}(U(x(t;x))) =P_t U(x)\leq e^{a t}U(x),
\ee
and  therefore $P_t$ is extendable to a linear operator on $\mathcal{F}_U$ (defined in \eqref{Eq_FU}) with norm $\|P_t\|\leq e^{a t}$.

The second inequality in \eqref{Eq_Lyap_cond} ensures that for any $t>t_0$,
\be\label{Eq_est0}
P_t U (x) \leq \kappa U(x) +b, \; \forall \; x\in \mathbb{R}^d.
\ee
By definition, a Markov semigroup is monotone, thus one may iterate  \eqref{Eq_est0} to obtain (by using $P_t \mathbbm{1}_{\mathbb{R}^d}=\mathbbm{1}_{\mathbb{R}^d}$), 
\be\label{Eq_est1}
P_{nt} U(x) \leq \kappa^n U(x) +\frac{b}{1-\kappa}, \; \; n\geq 1.
\ee
Consider now an arbitrary compact set $\mathfrak{B}$ in $\mathbb{R}^d$ and $f$ in $\mathcal{F}_{U}$, we have the bound
\bea
| P_{nt}f(x) -\mathbbm{1}_{\mathfrak{B}}(x) P_{nt} \mathbbm{1}_{\mathfrak{B}} f(x) | & \leq U(x) \underset{y\in \mathbb{R}^d\backslash \mathfrak{B}}\sup \frac{|P_{nt} f(y)|}{U(y)},\\
                                                     & \leq U(x) \| f\|_{U} \sup \frac{|P_{nt} U(y)|}{U(y)},
\eea
where we have used the basic inequality \eqref{Ineq_base} (with $p=1$). 
This last inequality  with \eqref{Eq_est1} leads to 
\be
| P_{nt}f(x) -\mathbbm{1}_{\mathfrak{B}}(x) P_{nt} \mathbbm{1}_{\mathfrak{B}} f(x) | \leq U(x) \| f\|_{U}  \left(\kappa^n + \frac{b}{1-\kappa}  \; \underset{y\in \mathcal{H}\backslash \mathfrak{B}}\sup \; \frac{1}{U(y)}\right).
\ee
Since $\lim_{|x|\rightarrow \infty} U(x)=\infty$, given $\epsilon>0$ and $n>1$ one may thus choose a compact set $\mathfrak{B}_n$ such that 
\be\label{Eq_compact_bounds}
\| P_{nt}f -\mathbbm{1}_{\mathfrak{B}_n}\, P_{nt} \mathbbm{1}_{\mathfrak{B}_n}f \|_U \leq \|f\|_{U} (\kappa +\epsilon)^n,
\ee
which leads to 
\be
\|P_{nt}-\mathbbm{1}_{\mathfrak{B}_n} \,P_{nt}\mathbbm{1}_{\mathfrak{B}_n}  \|_{\mathcal{L}(\mathcal{F}_U)} \leq (\kappa +\epsilon)^n.
\ee

{\bf Step 2.} We show now that the linear operator 
\be 
\Lambda =\mathbbm{1}_{\mathfrak{B}} \,P_{t} \mathbbm{1}_{\mathfrak{B}}: \mathcal{F}_U \longrightarrow \mathcal{F}_U,
\ee
 is  compact for any compact set $\mathfrak{B}$ of $\mathbb{R}^d$. This is equivalent to showing that for any sequence $g_k$ in 
$\mathcal{F}_U$ such that $\|g_k\|_U\leq 1$, one can extract a subsequence such that $\Lambda g_k$ is convergent in $\mathcal{F}_U$. 
Since $P_t$ is strongly Feller and $\mathbbm{1}_{\mathfrak{B}} g_k$ is bounded for each $k$, then $P_{t} \mathbbm{1}_{\mathfrak{B}} g_k$ belongs to  $\mathcal{C}_b(\mathfrak{B})$, by definition. Thus the sequence $(\Lambda g_k)$ lies in $\mathcal{C}(\mathfrak{B})$.

We have 
\be
|\Lambda g_k (x)|\leq \|g_k\|_U P_t U(x) \leq \kappa U(x)  + b \leq \kappa \; \underset{y\in \mathfrak{B}}\sup \; U(y) +b,\;\; x\in \mathfrak{B},
\ee
which shows that $\{\Lambda g_k\}$ is equibounded.  

Furthermore, since $P_t$ is strong Feller, it has a smooth kernel\footnote{A probability kernel  $ \mathfrak{T}_t$ allows for representing  the Markov semigroup $P_t$ as $P_t f(x)=\int \mathfrak{T}_t(x,\d y) f(y)$; e.g.~\cite[Prop.~1.2.3]{bakry2013analysis}. Having a smooth kernel means that $\mathfrak{T}_t(x,\d y)=\mathfrak{p}_t(x,y)\d y$ with $\mathfrak{p}_t$ infinitely differentiable, i.e.~smooth.} and we have for all $x$ and $x'$ in $\mathfrak{B}$
\bea
|\Lambda g_k (x)-\Lambda g_k (x')|&\leq \int_{y\in \mathfrak{B}}|\mathfrak{p}_t(x,y)- \mathfrak{p}_t(x',y)| |f(y)|\d y,\\
                                                         &\leq |x-x'| \underset{u,v\in \mathfrak{B}}\sup |\partial_u \mathfrak{p}_t(u,v)| \|g_k\|_U \int_{\mathfrak{B}} U(y) \d y,
\eea
which shows that $\{\Lambda g_k\}$ is equicontinuous. 

{\mk Thus, the Ascoli-Arzel\`a theorem \cite[p.~85]{yosida1995functional} applies and} guarantees that {\mk a subsequence from} $\Lambda g_k$ converges in $\mathcal{C}(\mathfrak{B})$ to $g$.  Now since $U\geq 1$,  {\mk the same extraction from} $\Lambda g_k$ converges to $g\mathbbm{1}_{\mathfrak{B}}$ in $\mathcal{F}_U$. We conclude that $\mathbbm{1}_{\mathfrak{B}} \,P_{t} \mathbbm{1}_{\mathfrak{B}}$ is a compact mapping for any compact set $\mathfrak{B}$ of $\mathbb{R}^d$.

{\bf Step 3.} Let $\mathfrak{B}_n$ be a sequence of compact sets satisfying \eqref{Eq_compact_bounds}, and let us consider the compact operators (from Step 2) $\mathfrak{C}_n$ defined by $\mathbbm{1}_{\mathfrak{B}_n} \,P_{n t} \mathbbm{1}_{\mathfrak{B}_n}$. We have then
\bea\label{Def_essnorm2}
\|P_{nt}\|_{ess}=\inf \Big\{ \|P_{nt}- \mathfrak{C}\|_{\mathcal{L}(\mathcal{F}_U)} \, : \, \mathfrak{C} \mbox{ is  a linear and compact operator of } \mathcal{F}_{U}\Big\}&  \leq \|P_{n t}- \mathfrak{C}_n\|\\
&\leq (\kappa +\epsilon)^n.
\eea
By applying to $P_t$ the Hadamard formula recalled in \eqref{Spec_rad}, we have  thus for $t>t_0$
\be
\mathbf{r}_{ess}(P_t)=\lim_{n\rightarrow \infty} \|P_{nt}\|_{ess}^{1/n}\leq \kappa +\epsilon,
\ee
for all $\epsilon>0,$ and we deduce \eqref{Eq_quasi}. 

The exponential convergence is then ensured by showing that there is no other eigenvalue
than 1 on the unit disk (or outside the unit disk) and that 1 is a simple eigenvalue; see \cite{bellet2006ergodic}.  \qed
\ep

\bibliographystyle{amsalpha}
\bibliography{atantet_rev_light}
\end{document}